# LIMIT THEOREMS FOR BIFURCATING MARKOV CHAINS. APPLICATION TO THE DETECTION OF CELLULAR AGING

By Julien Guyon

*ENPC-CERMICS*

We propose a general method to study dependent data in a binary tree, where an individual in one generation gives rise to two different offspring, one of type 0 and one of type 1, in the next generation. For any specific characteristic of these individuals, we assume that the characteristic is stochastic and depends on its ancestors' only through the mother's characteristic. The dependency structure may be described by a transition probability $P(x, dy\, dz)$ which gives the probability that the pair of daughters' characteristics is around $(y, z)$, given that the mother's characteristic is $x$. Note that $y$, the characteristic of the daughter of type 0, and $z$, that of the daughter of type 1, may be conditionally dependent given $x$, and their respective conditional distributions may differ. We then speak of bifurcating Markov chains.

We derive laws of large numbers and central limit theorems for such stochastic processes. We then apply these results to detect cellular aging in *Escherichia Coli*, using the data of Stewart et al. and a bifurcating autoregressive model.

## 1. Introduction.

1.1. *Motivation.* This study has been motivated by a collaboration [13] with biologists from the Laboratoire de Génétique Moléculaire, Évolutive et Médicale (INSERM U571, Faculté de Médecine Necker, Paris). F. Taddéi, E. J. Stewart, A. Lindner and G. Paul, together with R. Madden from the Institut des Hautes Études Scientifiques, have been working on *Escherichia Coli*'s aging. *E. Coli* is a single-celled, model organism. It has been widely studied by the biologists who have gathered a large amount of information on its physiology. Whereas aging is obvious in macroscopic organisms, it is not in single-celled ones, where, nevertheless, one has the best chances of describing and quantifying the molecular process involved. It is especially hard to identify in *E. Coli*, which reproduces without a juvenile phase and with an apparently symmetric division. Stewart et al. [22] have designed an experimental protocol which brings evidence of aging in *E. Coli* and we propose a statistical study of the data they collected.

In this section we describe the biological experiment and present the data (Section 1.2). Between-experiment averaging shows a clear segregation between the new- and old-pole derived progeny (see Section 1.2.1), whereas single-experiment data does not. In Section 1.3 we propose a linear Gaussian model that allows to study the populations of old and new poles experiment-wise. The model consists of a bifurcating Markovian dynamics. This motivates Section 2, where we give a detailed study of such stochastic processes. We pay special attention to limit theorems such as laws of large numbers (Theorems 8, 11 and 14) and central limit theorems (Theorem 19 and its corollaries). Eventually, in Section 3, we apply these results to the model, proving strong laws of large numbers and a central limit theorem (see Propositions 27 and 28), and derive rigorous estimation and test procedures which are performed on the data in order to detect cellular aging.

1.2. *The biology.* Here we briefly describe *E. Coli*'s life cycle, the experiment designed by Stewart et al. and the data they get. Figure 1 is taken from [22] where one can find further information.

1.2.1. *The experiment.* *E. Coli* is a rod-shaped bacterium. It reproduces by dividing in the middle, producing a new end per progeny cell (see Figure 1). This new end is called the new pole, whereas the other end is pre-existing and is called the old pole.

This defines an age in divisions for each pole, and hence, for each cell. One expects any cell component formed in the poles and with limited diffusion to accumulate at the old pole, so that there might be a physiological asymmetry between the old and new poles. To determine if *E. Coli* experiences aging related to the inheritance of the old pole, Stewart et al. followed 95 individual exponentially growing cells through up to nine generations in an automated fluorescence microscopy system which allowed them to determine the complete lineage, the identity of each pole and, among other physical parameters, the growth rate of each cell. Let us now present their results.

1.2.2. *Original data.* Each of the 95 films gives rise to a genealogical tree such as the ones in Figure 2. The new poles are the solid lines and the old poles the dashed lines. On the $y$-axis appears the growth rate, whereas the $x$-axis displays time in divisions. There is no striking evidence for reproductive asymmetry between the progeny cells visible to the naked eye. Note that generally the data is not regular: some generations are not completely observed, and in few cases a cell's growth rate might be measured whereas her sister's is not.

1.2.3. *Averaged data.* In order to eliminate the random effects which appear in Figure 2, Stewart et al. have averaged the 95 lineages by each unique cell position within the lineage. Figure 3 is the average tree thus produced for generations 5, 6 and 7. It clearly shows a segregation between the new and old poles. The old poles grow slower than the new poles, which is evidence for aging.

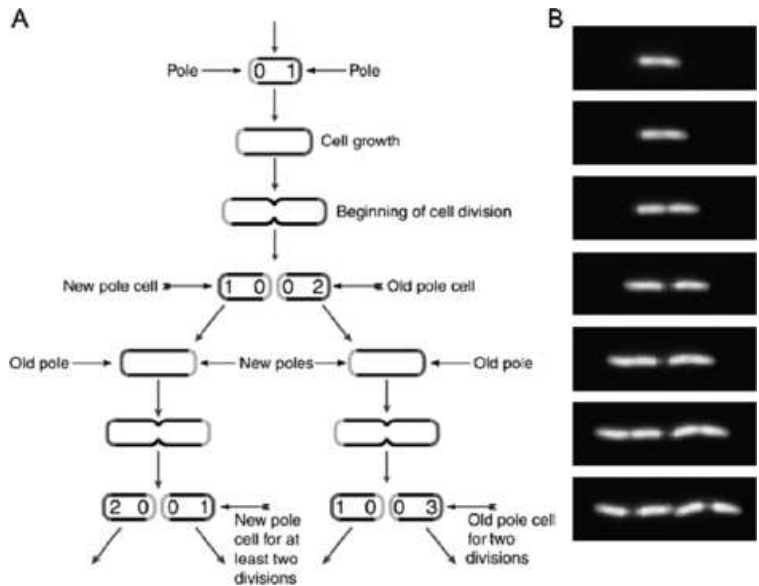

Fig. 1. *The life cycle of E. Coli, from [22].*

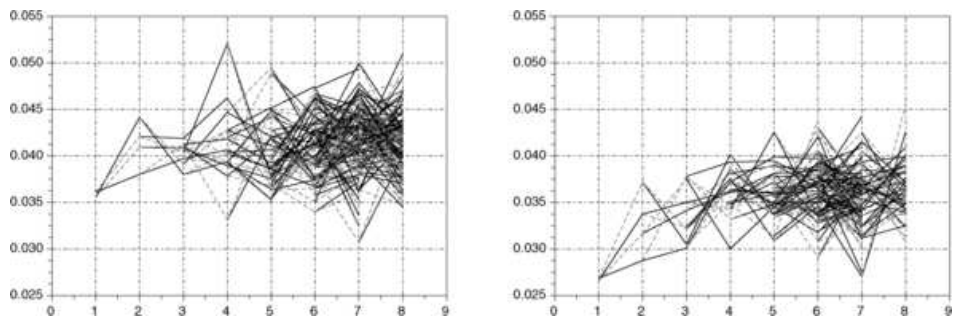

Fig. 2. *Two single-experiment data trees (two films).*

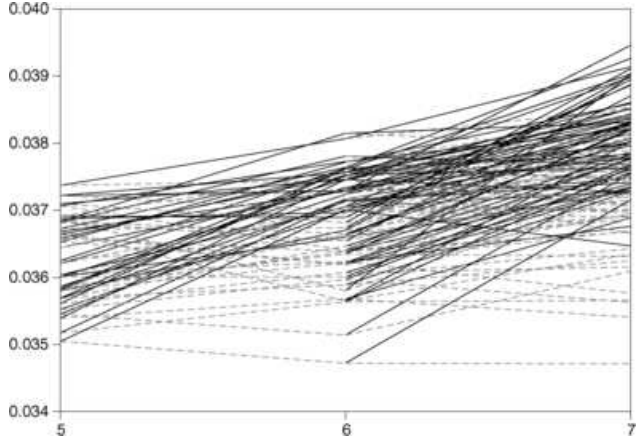

Fig. 3. *The average data tree*

However, we would prefer to study each experiment separately, since we do not know if the experiments are independent and/or identically distributed. Indeed, two initial cells giving rise to two different films are actually taken from the same macrocolony, so that there might be a correlation between the experiments. Moreover, as shown in Figure 2, the range of values of the growth rate changes from film to film, probably due to a slight change in the experimental conditions. In the next section we propose a statistical model that allows us to study the populations of old and new poles experiment-wise. It also has the advantage of taking into account the structure of the dependencies within a lineage. To be precise, contrary to Stewart et al., we take the effect of the mother into account, and we will prove that this effect is important (see Remark 39).

1.3. *The mathematical model.* In order to describe the dynamics of the growth rate, let $X_i$ denote the growth rate of individual $i$ and $n$ denote the mother of $2n$—the new pole progeny cell—and $2n+1$—the old pole progeny cell; see Figure 4. We propose the following Markovian model with memory one: $X_1$, the ancestor's growth rate, has distribution $\nu$ and for all $n \geq 1$,

$$
\begin{cases} X_{2n} = \alpha_0 X_n + \beta_0 + \varepsilon_{2n}, \\ X_{2n+1} = \alpha_1 X_n + \beta_1 + \varepsilon_{2n+1}, \end{cases} \quad (1)
$$

where $\alpha_0, \alpha_1 \in (-1, 1)$, $\beta_0, \beta_1 \in \mathbb{R}$ and $((\varepsilon_{2n}, \varepsilon_{2n+1}), n \geq 1)$ forms a sequence of independent and identically distributed (i.i.d.) centered bivariate Gaussian random variables (r.v.), say, $(\varepsilon_{2n}, \varepsilon_{2n+1}) \sim \mathcal{N}_2(0, \Gamma)$ with

$$
\Gamma = \sigma^2 \begin{pmatrix} 1 & \rho \\ \rho & 1 \end{pmatrix}, \qquad \sigma^2 > 0, \rho \in (-1, 1) \quad (2)
$$

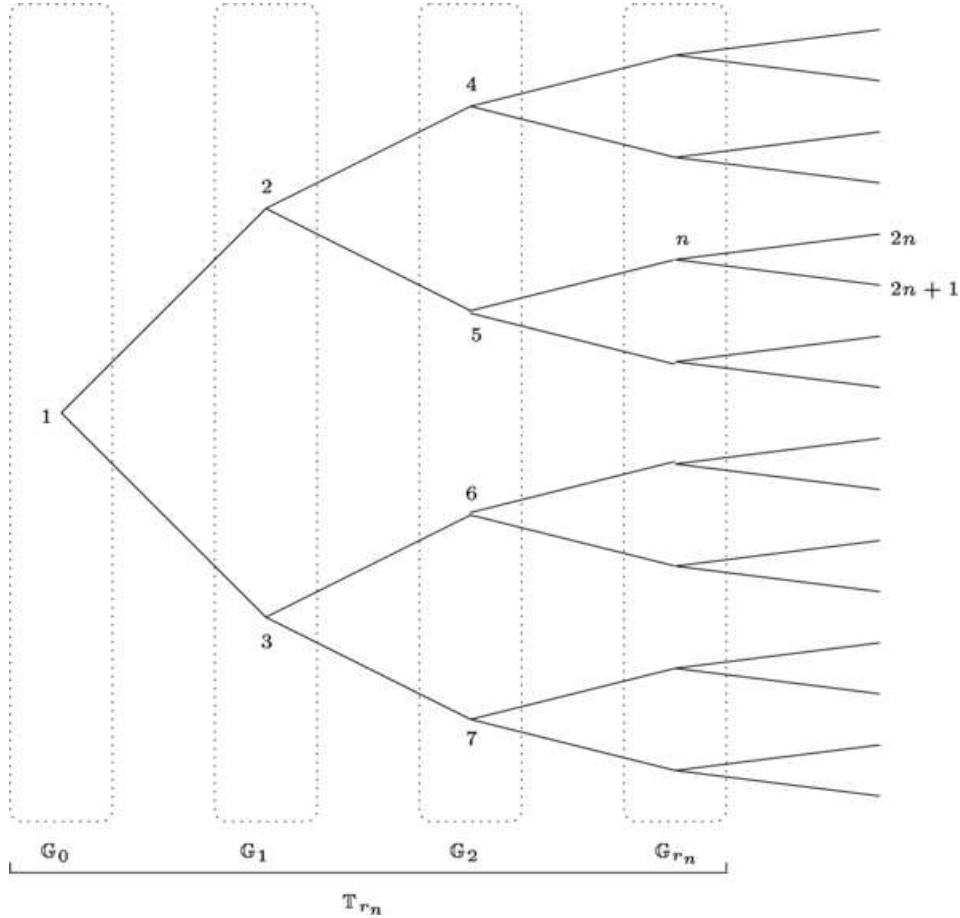

FIG. 4. *The binary tree* $\mathbb{T}$.

($\varepsilon_{2n}$ and $\varepsilon_{2n+1}$ are thus supposed to have common variance $\sigma^2$). We speak of memory one because a cell's growth rate is explained by its mother's. For instance, a Markovian model with memory two would also take into account the grandmother's growth rate. Considering Markovian models with memory two would allow us to test whether the grandmother effect is significant. In this article we concentrate on model (1) which we regard as "the simpler" reasonable model which describes a dependency within the colony.

REMARK 1. Since a Gaussian r.v. may take arbitrarily big negative values, here we allow the growth rate to take negative values. However, provided we correctly estimate the parameters, this should happen with extremely small probability.

We aim at the following:

(1) estimating the 4-dimensional parameter $\theta = (\alpha_0, \beta_0, \alpha_1, \beta_1)$, $\rho$ and $\sigma^2$,
(2) testing the null hypothesis $H_0 = \{(\alpha_0, \beta_0) = (\alpha_1, \beta_1)\}$ against its alternative $H_1 = \{(\alpha_0, \beta_0) \neq (\alpha_1, \beta_1)\}$.

In view of the biological question addressed here, point (2) is crucial: rejecting $H_0$ comes down to accepting that the dynamics of the growth rate of the old pole offspring is different from that of the new pole offspring. We shall actually see that the old pole progeny cell experiences slowed growth rate and, hence, should be considered an aging parent repeatedly producing rejuvenated offspring.

Bifurcating autoregressive (BAR) models, such as model (1), have already been studied. Cowan and Staudte [9] were pioneers and studied model (1) in the special case when $(\alpha_0, \beta_0) = (\alpha_1, \beta_1)$, that is, under $H_0$. Several extensions [2, 3, 4, 8, 14, 15, 16, 21] followed, improving inference results or/and generalizing the model, but no distinction was ever made between new and old poles. In mathematical terms, in all these articles, the distribution of $X_{2n}$ given $X_n$ has always been assumed to be the same as the distribution of $X_{2n+1}$ given $X_n$. Now, detecting a discrepancy between these distributions is the central question addressed here. Hence, model (1) generalizes existing BAR models and allows us to detect dissymmetry between sisters. Such a generalization is a source of mathematical difficulties. For instance, there is no stationary distribution in the sense of [9], that is, to say a distribution common to all cells in the clone. Besides, $(X_n, n \geq 1)$ does not converge in distribution. For example, by computing characteristic functions, it is easy to see that the sequence of always new poles' growth rates $(X_{2^n}, n \geq 1)$ and that of always old poles' growth rates $(X_{2^{n+1}-1}, n \geq 1)$ both converge in distribution, but, unless $\alpha_0^2 = \alpha_1^2$ and $\beta_0/(1-\alpha_0) = \beta_1/(1-\alpha_1)$, the (Gaussian) limit distributions are distinct. This leads us to develop a new theory (see Section 2). Note that in a recent work Evans and Steinsaltz [11] also address the question of dissymmetry between sisters' growth rates, by proposing a superprocess model for damage segregation and showing that optimal growth is achieved by unequal division of (deterministically accumulating) damage between the daughters.

We shall call $X = (X_n, n \geq 1)$ a bifurcating Markov chain. The next section is devoted to the study of this family of stochastic processes. Establishing laws of large numbers and central limit theorems will be crucial in achieving the two above objectives. That is the reason why we will pay special attention to such limit theorems.

**2. Bifurcating Markov chains. Limit theorems.**

2.1. *Definitions.* Markov chains (MCs) are usually indexed by the integers. Here we give a definition of a MC when the index set is the (regular) binary tree $\mathbb{T}$. We then speak of a bifurcating Markov chain or a $\mathbb{T}$-Markov chain, which we often write $\mathbb{T}$-MC. $\mathbb{T}$-MCs are well adapted to modeling data on the descent of an initial individual, where each individual in one generation gives rise to two offspring in the next one. Cell lineage data, such as the one presented in Section 1.2, are typically of this kind.

Let us introduce some notation about the binary tree $\mathbb{T}$; see Figure 4. Each vertex $n \in \mathbb{T}$ is seen as a positive integer $n \in \mathbb{N}^*$. It should be thought of as an individual or a cell. It has exactly two daughters, $2n$ and $2n+1$, and we label the root 1. We denote by $\mathbb{G}_q = \{2^q, 2^q + 1, \ldots, 2^{q+1} - 1\}$ the $q$th generation and by $\mathbb{T}_r = \bigcup_{q=0}^r \mathbb{G}_q$ the subtree consisting of the first $r+1$ generations. With this notation, $\mathbb{G}_0 = \{1\}$ and, $|\cdot|$ standing for the cardinality, $|\mathbb{G}_q| = 2^q$ and $|\mathbb{T}_r| = 2^{r+1} - 1$. We also denote by $r_n = \lfloor \log_2 n \rfloor$ the generation of individual $n$, that is, $n \in \mathbb{G}_{r_n}$. In terms of labeling the vertices, $\mathbb{T}$ is assimilated to $\mathbb{N}^*$, but the topology is different: within $\mathbb{T}$, $n$ and $2n$ (resp. $n$ and $2n+1$) should be seen as neighbors.

Let $(S, \mathcal{S})$ be a metric space endowed with its Borel $\sigma$-field, and think of it as the state space. For instance, in the BAR model (1), $S = \mathbb{R}$. For any integer $p \geq 2$, we equip $S^p$ with the product $\sigma$-field, say, $\mathcal{S}^p$.

DEFINITION 2. We call $\mathbb{T}$-transition probability any mapping $P : S \times \mathcal{S}^2 \to [0, 1]$ such that:

- $P(\cdot, A)$ is measurable for all $A \in \mathcal{S}^2$,
- $P(x, \cdot)$ is a probability measure on $(S^2, \mathcal{S}^2)$ for all $x \in S$.

We also define, for $x \in S$ and $B \in \mathcal{S}$, $P_0(x, B) = P(x, B \times S)$ and $P_1(x, B) = P(x, S \times B)$. $P_0$ and $P_1$ are transition probabilities on $(S, \mathcal{S})$. In the BAR model (1), they respectively correspond to the transition probabilities of the new poles and of the old poles.

For $p \geq 1$, we denote by $\mathcal{B}(S^p)$ [resp. $\mathcal{B}_b(S^p)$, $\mathcal{C}(S^p)$, $\mathcal{C}_b(S^p)$] the set of all $\mathcal{S}^p$-measurable (resp. $\mathcal{S}^p$-measurable and bounded, continuous, continuous and bounded) mappings $f : S^p \to \mathbb{R}$. If $p \in \{2, 3\}$ and $f \in \mathcal{B}(S^p)$, when it is defined, we denote by $Pf \in \mathcal{B}(S)$ the function

$$x \mapsto Pf(x) = \begin{cases} \int_{S^2} f(y, z) P(x, dy\, dz), & \text{if } p = 2, \\ \int_{S^2} f(x, y, z) P(x, dy\, dz), & \text{if } p = 3. \end{cases}$$

Let $(\Omega, \mathcal{F}, (\mathcal{F}_r, r \in \mathbb{N}), \mathbb{P})$ be a filtered probability space and, defined on it, $(X_n, n \in \mathbb{T})$ be a family of $S$-valued random variables. Let $\nu$ be a probability on $(S, \mathcal{S})$ and $P$ be a $\mathbb{T}$-transition probability.

DEFINITION 3. We say that $(X_n, n \in \mathbb{T})$ is a $(\mathcal{F}_r)$-bifurcating Markov chain, or $(\mathcal{F}_r)$-$\mathbb{T}$-MC (with initial distribution $\nu$ and $\mathbb{T}$-transition probability $P$), if:

- $X_n$ is $\mathcal{F}_{r_n}$-measurable for all $n \in \mathbb{T}$,
- $X_1$ has distribution $\nu$,
- for all $q \in \mathbb{N}$ and for all family $(f_n, n \in \mathbb{G}_q)$ in $\mathcal{B}_b(S^2)$,

$$\mathbb{E}\left[\prod_{n \in \mathbb{G}_q} f_n(X_{2n}, X_{2n+1}) \Big| \mathcal{F}_q\right] = \prod_{n \in \mathbb{G}_q} Pf_n(X_n).$$

This means that, given generations 0 to $q$, $\mathbb{T}_q$, one builds generation $\mathbb{G}_{q+1}$ by drawing $2^q$ *independent* couples $(X_{2n}, X_{2n+1})$ according to $P(X_n, \cdot)$ ($n \in \mathbb{G}_q$). A $(\mathcal{F}_r)$-$\mathbb{T}$-MC is also a $(\mathcal{F}_r^X)$-$\mathbb{T}$-MC, where $\mathcal{F}_r^X = \sigma(X_i, i \in \mathbb{T}_r)$. When unstated, the filtration implicitly used is the latter. Note that for $f \in \mathcal{B}_b(S^2)$, $\mathbb{E}[f(X_{2n}, X_{2n+1})|\mathcal{F}_{r_n}]$ factorizes through the random variable $X_n$, so that $\mathbb{E}[f(X_{2n}, X_{2n+1})|\mathcal{F}_{r_n}] = \mathbb{E}[f(X_{2n}, X_{2n+1})|X_n]$. This means that any $X_n$ depends on past generations only through her mother. This explains why we speak of a Markov chain (with memory one).

Last but not least, note that, contrary to much of the (still sparse) literature on the subject, we allow conditional dependency between sisters. Conditional independence corresponds to the case when $P$ factorizes as a product $P_0 \otimes P_1$ of two transition probabilities on $(S, \mathcal{S})$, that is, $P(x, dy\,dz) = P_0(x, dy) \otimes P_1(x, dz)$ for all $x \in S$. [1, 5, 17, 18, 19, 23, 24] deal with more general than binary—and possibly random—trees, but all assume that, conditionally on their mother's type, the daughters have independent and identically distributed types. In our case, this corresponds to conditional independency with $P_0 = P_1$. As said in Section 1.3, to our best knowledge, BAR models, although they allow for conditional dependence, have always been studied until now under the assumption that $P_0 = P_1$. Now, detecting that $P_0 \neq P_1$ will be the central question when we study cellular aging (see Section 3). Moreover, like in [1], we consider general state spaces whereas [5] deal only with countable ones and [17, 18, 19, 23, 24] only with finite ones. Note that in the latter case, one may regard a $\mathbb{T}$-MC $X$ as a multitype branching process and apply Maâouia and Touati's identification techniques [19].

2.2. *Weak law of large numbers.*

2.2.1. *Introduction.* A first natural question is to know whether a $\mathbb{T}$-Markov chain $(X_n)$ obeys laws of large numbers (LLN), that is, convergence of empirical means.

Given $f \in \mathcal{B}(S)$ and a finite subset $I \subset \mathbb{T}$, let us write $M_I(f) = \sum_{i \in I} f(X_i)$ and $\overline{M}_I(f) = |I|^{-1} M_I(f)$. Several empirical averages can be considered:

- One may average over the $q$th generation, that is, compute $\overline{M}_{\mathbb{G}_q}(f)$.
- One may prefer to average over the first $r+1$ generations, that is, compute $\overline{M}_{\mathbb{T}_r}(f)$. This is meaningful because $\mathbb{G}_q$ naturally precedes $\mathbb{G}_{q+1}$, since one cannot draw the whole $(q+1)$th generation without having completely drawn the $q$th one.
- One may also average over the "first" $n$ individuals, that is, compute $n^{-1}\sum_{i=1}^{n} f(X_i)$. However, there is no natural order within a generation $\mathbb{G}_q$: all the individuals $(X_n, n \in \mathbb{G}_q)$ of the $q$th level can be generated simultaneously. That is why we introduce the set $\mathfrak{S}$ of all permutations of $\mathbb{N}^*$ that leave each $\mathbb{G}_q$ invariant and, for $f \in \mathcal{B}(S)$ and $\pi \in \mathfrak{S}$, consider the sums

$$M_n^\pi(f) = \sum_{i=1}^{n} f(X_{\pi(i)}).$$

As far as the asymptotic behavior of $\overline{M}_n^\pi(f) \equiv n^{-1} M_n^\pi(f)$ is concerned, the choice of $\pi$ matters. This is because each new generation is essentially the same size as the total of all previous ancestors. To illustrate this, consider the following example.

EXAMPLE 4. Assume that $S = \{0,1\}$, $f = \mathrm{id}_S$ and, whatever the mother's type, $X_{2n} = 1$ and $X_{2n+1} = 0$—in other words, $P(x, dy\, dz) = \delta_1(dy)\delta_0(dz)$ for all $x \in S$, where $\delta_x$ stands for the Dirac mass at point $x$. If $\pi \in \mathfrak{S}$ sends the "first half" of each $\mathbb{G}_q$, that is, $\{2^q, 2^q+1, \ldots, 3 \cdot 2^{q-1}-1\}$, onto the even elements in $\mathbb{G}_q$, that is, $\mathbb{G}_q \cap (2\mathbb{N})$, then $\liminf_{n\to\infty} \overline{M}_n^\pi(f) = 1/2$ and $\limsup_{n\to\infty} \overline{M}_n^\pi(f) = 2/3$: the empirical average $\overline{M}_n^\pi(f)$ oscillates between $1/2$ and $2/3$. Conversely, if $\pi$ sends the "first half" of each $\mathbb{G}_q$ onto the odd elements in $\mathbb{G}_q$, that is, $\mathbb{G}_q \cap (2\mathbb{N}+1)$, then $\liminf_{n\to\infty} \overline{M}_n^\pi(f) = 1/3$ and $\limsup_{n\to\infty} \overline{M}_n^\pi(f) = 1/2$. But for $\pi = \mathrm{id}_{\mathbb{N}^*}$, $\overline{M}_n^\pi(f)$ converges to $1/2$.

- A natural answer to this issue is to explore each new generation "by chance," that is, to draw a permutation $\Pi$ "uniformly" on $\mathfrak{S}$, independently on $X = (X_n, n \in \mathbb{T})$. Drawing $\Pi$ "uniformly" on $\mathfrak{S}$ means drawing the restriction of $\Pi$ on $\mathbb{G}_q$ uniformly among the $(2^q)!$ permutations of $\mathbb{G}_q$, independently for each $q$. Then we consider the empirical average $\overline{M}_n^\Pi(f) \equiv n^{-1} M_n^\Pi(f)$, where

$$M_n^\Pi(f) = \sum_{i=1}^{n} f(X_{\Pi(i)}).$$

Thus, we introduce extra randomness, but this will allow us to get through Liapunov's condition when we try to derive a central limit theorem for $X$.

REMARK 5. Note that for all $\pi \in \mathfrak{S}$ and $r \geq 0$, $M_{|\mathbb{T}_r|}^{\pi}(f) = M_{\mathbb{T}_r}(f)$. Besides, for all $r \geq 0$,

$$\overline{M}_{\mathbb{T}_r}(f) = \sum_{q=0}^{r} \frac{|\mathbb{G}_q|}{|\mathbb{T}_r|} \overline{M}_{\mathbb{G}_q}(f) \tag{3}$$

and, for all $\pi \in \mathfrak{S}$ and $n \geq 1$,

$$\frac{1}{n} M_n^{\Pi}(f) = \sum_{q=0}^{r_n - 1} \frac{|\mathbb{G}_q|}{n} \overline{M}_{\mathbb{G}_q}(f) + \frac{1}{n} \sum_{i=2^{r_n}}^{n} f(X_{\Pi(i)}) \tag{4}$$

(we systematically use the convention that a sum over an empty set is zero).

Because of the branching, empirical averages of $\mathbb{T}$-MCs may not behave like corresponding MCs ones. Precisely, given a transition probability $R$, a LLN may hold for MCs with transition probability $R$ and fail for $\mathbb{T}$-MCs with $\mathbb{T}$-transition probability $R \otimes R$. A very simple but crucial illustration of this is Example 6. At least for the case of finite state spaces, let us keep in mind that periodicity is problematic and that there is more to ask than irreducibility and recurrence for a $\mathbb{T}$-MC to obey a LLN.

EXAMPLE 6. Consider the two-state MC, say, $S = \{0, 1\}$, with $R(0, \cdot) = \delta_1$ and $R(1, \cdot) = \delta_0$. For a MC $Y$, a LLN holds true: whatever the initial state, the empirical average $\frac{1}{n} \sum_{i=1}^{n} Y_i$ converges to $1/2$ when $n$ grows to infinity. But for the corresponding $\mathbb{T}$-MC $X$ with $\mathbb{T}$-transition probability $P = R \otimes R$, that is, $P(0, \cdot) = \delta_1 \otimes \delta_1$ and $P(1, \cdot) = \delta_0 \otimes \delta_0$, it endlessly fluctuates from $1/3$ to $2/3$.

2.2.2. *Results.* Here we ask ourselves whether the various empirical averages introduced in Section 2.2.1 converge, in quadratic mean, when the size of the tree grows to infinity. We then speak of weak LLNs. A sufficient condition for a weak LLN to hold appears to be the ergodicity (see Definition 7) of the induced MC $(Y_r, r \in \mathbb{N})$, defined as follows. Start from the root and recursively choose one of the two daughters tossing a balanced coin, independently on the $\mathbb{T}$-MC $X$. In mathematical terms, $Y_0 = X_1$ and if $Y_r = X_n$, then $Y_{r+1} = X_{2n + \zeta_{r+1}}$ for a sequence of independent balanced Bernoulli r.v. $(\zeta_q, q \in \mathbb{N}^*)$ independent on $(X, \Pi)$. Here, "balanced" means that $\mathbb{P}(\zeta_q = 0) = \mathbb{P}(\zeta_q = 1) = 1/2$. It is easy to check that $(Y_r, r \in \mathbb{N})$ is a MC with initial distribution $\nu$ and transition probability

$$Q = \frac{P_0 + P_1}{2}.$$

DEFINITION 7. We say that a MC $Y$ is ergodic if there exists a probability $\mu$ on $(S, \mathcal{S})$ such that $\lim_{r \to \infty} \mathbb{E}_x[f(Y_r)] = \int_S f \, d\mu$ for all $x \in S$ and $f \in \mathcal{C}_b(S)$.

Then $\mu$ is the unique stationary distribution of $Y$, and the sequence $(Y_r, r \in \mathbb{N})$ converges in distribution to $\mu$. Sufficient conditions for ergodicity may be found in [6, 20]. We are now in the position to state the main theorem of this section:

THEOREM 8. *Assume that the induced MC $(Y_r, r \in \mathbb{N})$ is ergodic, with stationary distribution $\mu$. Then, for any $f \in \mathcal{C}_b(S)$, the three empirical averages $\overline{M}_{\mathbb{G}_q}(f)$, $\overline{M}_{\mathbb{T}_r}(f)$ and $\overline{M}_n^{\Pi}(f)$ converge to $(\mu, f)$ in $L^2$.*

REMARK 9. It is noteworthy that the asymptotic behavior of the three above empirical averages depends on the $\mathbb{T}$-transition probability $P$ only through $Q = (P_0 + P_1)/2$.

REMARK 10. Athreya and Kang [1] use an analogous ergodicity hypothesis to get laws of large numbers. Namely, their results hold for Galton–Watson trees in which particles move according to a Markov chain $R$ on $(S, \mathcal{S})$, and they assume $\lim_{m \to \infty} R^m(x, \cdot) = \mu$. If this happens uniformly in $x$ on the compact subsets of $S$, they get a weak LLN; if this happens uniformly in $x$ on $S$, they get a strong LLN. Observe that we do not assume any uniformity in $x$ (but our tree is deterministic).

In the application (Section 3) the function $f$ will typically be unbounded so that we shall actually prove an extended version of Theorem 8. To this end, let us first introduce some notation. We denote by:

- $i \wedge j$ the most recent common ancestor of $i, j \in \mathbb{T}$,
- $f \otimes g$ the mapping $(x, y) \mapsto f(x)g(y)$,
- $Q^p$ the $p$th iterated of $Q$, recursively defined by the formulas $Q^0(x, \cdot) = \delta_x$ and $Q^{p+1}(x, B) = \int_S Q(x, dy) Q^p(y, B)$ for all $B \in \mathcal{S}$; $Q^p$ is a transition probability on $(S, \mathcal{S})$,
- $\nu Q$ the distribution on $(S, \mathcal{S})$ defined by $\nu Q(B) = \int_S \nu(dx) Q(x, B)$; $\nu Q^p$ is the law of $Y_p$,
- $(Qf)(x) = \int_S f(y) Q(x, dy)$ when it is defined,
- $\nu(f)$ or $(\nu, f)$ the integral $\int_S f \, d\nu$ when it is defined.

With such a notation, for any distribution $\lambda$, transition probability $Q$ and function $f \in \mathcal{B}(S)$ such that $\lambda Q(|f|) < \infty$, we have $\lambda Q(f) = \lambda(Qf)$ which is, hence, simply written $\lambda Q f$.

Now, let $F$ denote a subspace of $\mathcal{B}(S)$ such that:

(i) $F$ contains the constants,
(ii) $F^2 \subset F$,
(iii) $F \otimes F \subset L^1(P(x, \cdot))$ for all $x \in S$, and $P(F \otimes F) \subset F$,

(iv) there exists a probability measure $\mu$ on $(S, \mathcal{S})$ such that $F \subset L^1(\mu)$ and $\lim_{r \to \infty} \mathbb{E}_x[f(Y_r)] = (\mu, f)$ for all $x \in S$ and $f \in F$,

(v) for all $f \in F$, there exists $g \in F$ such that for all $r \in \mathbb{N}$, $|Q^r f| \leq g$,

(vi) $F \subset L^1(\nu)$,

where we have used the notation $F^2 = \{f^2 | f \in F\}$, $F \otimes F = \{f \otimes g | f, g \in F\}$ and $PE = \{Pf | f \in E\}$ whenever an operator $P$ acts on a set $E$. Note that (i) and (iii) imply the condition

(iii′) for $z = 0, 1$, $F \subset L^1(P_z(x, \cdot))$ for all $x \in S$ and $P_z F \subset F$,

since $P_0 f = P(f \otimes \mathbf{1})$ and $P_1 f = P(\mathbf{1} \otimes f)$. This in its turn implies

(iii″) $F \subset L^1(Q(x, \cdot))$ for all $x \in S$ and $QF \subset F$,

so that in (iv) and (v) $\mathbb{E}_x[f(Y_r)] = Q^r f(x)$ is well defined. Note also that if $F$ contains enough functions, that is, if it contains the set of all bounded Lipschitz functions, then $\mu$ is the unique stationary distribution of $Y$, that is, $\mu Q = \mu$. The next theorem states that the result in Theorem 8 remains true for $f$'s in such a $F$:

THEOREM 11. *Let $F$ satisfy conditions* (i)–(vi) *above. Then, for any $f \in F$, the three empirical averages $\overline{M}_{\mathbb{G}_q}(f)$, $\overline{M}_{\mathbb{T}_r}(f)$ and $\overline{M}_n^\Pi(f)$ converge to $(\mu, f)$ in $L^2$.*

Obviously $F = \mathcal{C}_b(S)$ fulfills conditions (i)–(vi) as soon as $Y$ is ergodic, so that Theorem 11 implies Theorem 8. In Section 3 we shall take $F$ to be the set of all continuous and polynomially growing functions.

We shall also need an easy extension of Theorem 11 to the case when $f$ does not only depend on an individual $X_i$, but on the mother–daughters triangle $(X_i, X_{2i}, X_{2i+1})$. This will be useful in the application (Section 3). Let us denote $\Delta_n = (X_n, X_{2n}, X_{2n+1})$ and, for $f \in \mathcal{B}(S^3)$ and $I$ a finite subset of $\mathbb{T}$,

$$M_I(f) = \sum_{i \in I} f(\Delta_i) \quad \text{and} \quad M_n^\Pi(f) = \sum_{i=1}^n f(\Delta_{\Pi(i)}).$$

Then we have the following:

THEOREM 12. *Let $F$ satisfy conditions* (i)–(vi). *Let $f \in \mathcal{B}(S^3)$ such that $Pf$ and $Pf^2$ exist and belong to $F$. Then the three empirical averages $\overline{M}_{\mathbb{G}_q}(f)$, $\overline{M}_{\mathbb{T}_r}(f)$ and $\overline{M}_n^\Pi(f)$ converge to $(\mu, Pf)$ in $L^2$.*

2.2.3. *Proofs.* This section is devoted to the proofs of Theorems 11 and 12.

PROOF OF THEOREM 11. Considering the function $f - (\mu, f)$ leaves us with the case when $(\mu, f) = 0$. Then condition (iv) implies that

$$\forall x \in S \quad \lim_{r \to \infty} Q^r f(x) = 0. \tag{5}$$

We shall study the three empirical averages $\overline{M}_{\mathbb{G}_q}(f)$, $\overline{M}_{\mathbb{T}_r}(f)$ and $\overline{M}_n^\Pi(f)$ successively.

*Step* 1. Let us first deal with $\overline{M}_{\mathbb{G}_q}(f)$. First note that $f(X_i) \in L^2$ for all $i \in \mathbb{G}_q$. Indeed, there is a unique path $(z_1, \ldots, z_q) \in \{0,1\}^q$ in the binary tree from the root 1 to $i$; here $(z_1, \ldots, z_q)$ should be seen as the realization of the coin toss r.v. $(\zeta_1, \ldots, \zeta_q)$ that joins 1 to $i$. For instance, $(1,0,0,1)$ is the path from 1 to 25. Thus,

$$\mathbb{E}[f(X_i)^2] = \nu P_{z_1} \cdots P_{z_q} f^2,$$

which, from (ii), (iii') and (vi), is finite.

Independently on $X$, let us now draw two independent indices $I_q$ and $J_q$ uniformly from $\mathbb{G}_q$. Then $f(X_{I_q})f(X_{J_q}) \in L^1$ and we have

$$\mathbb{E}[\overline{M}_{\mathbb{G}_q}(f)^2] = \frac{1}{|\mathbb{G}_q|^2} \sum_{i,j \in \mathbb{G}_q} \mathbb{E}[f(X_i)f(X_j)]$$
$$= \mathbb{E}[f(X_{I_q})f(X_{J_q})].$$

Let us fix $p \in \{0, \ldots, q\}$ and reason conditionally on the event $\{I_q \wedge J_q \in \mathbb{G}_p\}$. Then $I_q \wedge J_q$ is uniformly distributed on $\mathbb{G}_p$, so that $X_{I_q \wedge J_q}$ has the same distribution as $Y_p$, that is, has distribution $\nu Q^p$. Besides, for $p < q$, conditionally on the states $(X_{2(I_q \wedge J_q)}, X_{2(I_q \wedge J_q)+1})$ of the two daughters of $I_q \wedge J_q$, $X_{I_q}$ and $X_{J_q}$ are independent and have the same distribution as $Y_{q-p-1}$ with respective initial conditions $X_{2(I_q \wedge J_q)}$ and $X_{2(I_q \wedge J_q)+1}$. Provided we use the convention that $P(Q^{-1}f \otimes Q^{-1}f) = f^2$, we then have

$$\mathbb{E}[f(X_{I_q})f(X_{J_q})|I_q \wedge J_q \in \mathbb{G}_p] = \nu Q^p P(Q^{q-p-1}f \otimes Q^{q-p-1}f). \tag{6}$$

Now, $\mathbb{P}(I_q \wedge J_q \in \mathbb{G}_q) = \mathbb{P}(I_q = J_q) = 2^{-q}$ and, for $p \in \{0, \ldots, q-1\}$, $\mathbb{P}(I_q \wedge J_q \in \mathbb{G}_p) = 2^{-p-1}$. Indeed, since $I_q$ and $J_q$ are independent, the paths $(\zeta_1^I, \ldots, \zeta_q^I)$ from 1 to $I_q$ and $(\zeta_1^J, \ldots, \zeta_q^J)$ from 1 to $J_q$ are independent so that for $p \in \{0, \ldots, q-1\}$,

$$\mathbb{P}(I_q \wedge J_q \in \mathbb{G}_p) = \mathbb{P}(\zeta_1^I = \zeta_1^J, \ldots, \zeta_p^I = \zeta_p^J, \zeta_{p+1}^I \neq \zeta_{p+1}^J) = 2^{-p-1}.$$

In short, we write $\mathbb{P}(I_q \wedge J_q \in \mathbb{G}_p) = 2^{-p-\mathbf{1}_{\{p<q\}}}$. Combined with (6), this finally gives

$$
\begin{aligned}
(7) \quad \mathbb{E}[\overline{M}_{\mathbb{G}_q}(f)^2] &= \mathbb{E}[f(X_{I_q})f(X_{J_q})] \\
&= \sum_{p=0}^{q} 2^{-p-\mathbf{1}_{\{p<q\}}} \nu Q^p P(Q^{q-p-1}f \otimes Q^{q-p-1}f).
\end{aligned}
$$

Let us now fix $\varepsilon > 0$ and choose $p_\varepsilon \in \mathbb{N}$ such that $2^{-p_\varepsilon} \leq \varepsilon$. Then $\sum_{p>p_\varepsilon} 2^{-p} \leq \varepsilon$. Besides, from (iii), (v) and (vi), there is a $c_f \geq 0$ such that

$$
(8) \qquad |\nu Q^p P(Q^{q-p-1}f \otimes Q^{q-p-1}f)| \leq c_f
$$

for all $0 \leq p \leq q$. Hence, for $q > p_\varepsilon$,

$$
(9) \qquad \mathbb{E}[\overline{M}_{\mathbb{G}_q}(f)^2] \leq \varepsilon c_f + \sum_{p=0}^{p_\varepsilon} |\nu Q^p P(Q^{q-p-1}f \otimes Q^{q-p-1}f)|.
$$

From (v), there exists $g \in F$ such that, for all $r \in \mathbb{N}$, $|Q^r f \otimes Q^r f| \leq g \otimes g$. From (iii), $P(g \otimes g) \in F$ so that, using (v) and (vi), $P(g \otimes g) \in \bigcap_{p \in \mathbb{N}} L^1(\nu Q^p)$. This shows that the sequence of functions $(Q^r f \otimes Q^r f, r \in \mathbb{N})$ is dominated by $g \otimes g \in \bigcap_{p \in \mathbb{N}} L^1(\nu Q^p P)$. Then (5) and Lebesgue's dominated convergence theorem imply that

$$
(10) \qquad \forall p \in \mathbb{N} \qquad \lim_{r \to \infty} \nu Q^p P(Q^r f \otimes Q^r f) = 0.
$$

As a consequence, the r.h.s. of (9) converges to $\varepsilon c_f$ as $q$ grows to infinity. Since $\varepsilon$ is arbitrary, the proof is complete for $\overline{M}_{\mathbb{G}_q}(f)$.

Convergence results for $\overline{M}_{\mathbb{T}_r}(f)$ or $\overline{M}_n^\Pi(f)$ may be easily deduced from those for $\overline{M}_{\mathbb{G}_q}(f)$ by using (3) or (4) and the following lemma.

LEMMA 13. *Let $(u_r, r \in \mathbb{N})$ be a sequence of nonnegative real numbers converging to 0. Let*

$$
v_r = \sum_{q=0}^{r} \frac{|\mathbb{G}_q|}{|\mathbb{T}_r|} u_q \quad \text{and} \quad a_n = \sum_{q=0}^{r_n - 1} \frac{|\mathbb{G}_q|}{n} u_q.
$$

*Then $(v_r, r \in \mathbb{N})$ and $(a_n, n \in \mathbb{N}^*)$ converge to 0.*

PROOF. Let us fix $\varepsilon > 0$. We can find $q_\varepsilon \in \mathbb{N}^*$ such that $u_q \leq \varepsilon$ for all $q \geq q_\varepsilon$. Letting $M = \sup_{q \in \mathbb{N}} u_q$, we then have, for all $r \geq q_\varepsilon$, $v_r \leq \varepsilon + M \sum_{q=0}^{q_\varepsilon - 1} \frac{|\mathbb{G}_q|}{|\mathbb{T}_r|}$. The r.h.s. tends to $\varepsilon$ as $r$ grows to infinity, so that $\lim_{r \to \infty} v_r = 0$.

As for $(a_n, n \in \mathbb{N}^*)$, it is enough to notice that, since $|\mathbb{T}_{r_n - 1}| \leq n$, $a_n = v_{r_n - 1} |\mathbb{T}_{r_n - 1}|/n \leq v_{r_n - 1}$ and to apply the result for $(v_r, r \in \mathbb{N})$. $\square$

*Step* 2. Let us now treat $\overline{M}_{\mathbb{T}_r}(f)$. From (3), we have by the triangle inequality $\|\overline{M}_{\mathbb{T}_r}(f)\|_{L^2} \leq \sum_{q=0}^{r} \frac{|\mathbb{G}_q|}{|\mathbb{T}_r|} \|\overline{M}_{\mathbb{G}_q}(f)\|_{L^2}$. From Step 1, $\overline{M}_{\mathbb{G}_q}(f)$ converges to 0 in quadratic mean. Lemma 13 implies that the r.h.s. tends to 0, which ends the proof for $\overline{M}_{\mathbb{T}_r}(f)$.

*Step* 3. Eventually, let us look at $\overline{M}_n^\Pi(f)$. From (4) and the triangle inequality, $\|\overline{M}_n^\Pi(f)\|_{L^2} \leq a_n + b_n$, where

$$a_n = \sum_{q=0}^{r_n-1} \frac{|\mathbb{G}_q|}{n} \|\overline{M}_{\mathbb{G}_q}(f)\|_{L^2} \quad \text{and} \quad b_n = \left\| \frac{1}{n} \sum_{i=2^{r_n}}^{n} f(X_{\Pi(i)}) \right\|_{L^2}.$$

Since $\overline{M}_{\mathbb{G}_q}(f)$ converges to 0 in quadratic mean, Lemma 13 implies that $\lim_{n\to\infty} a_n = 0$. As for $b_n$, first note that since each $f(X_i)$ belongs to $L^2$, $f(X_{\Pi(i)})f(X_{\Pi(j)}) \in L^1$ for all $i, j \in \{2^{r_n}, \ldots, n\}$ so that

$$b_n^2 = \frac{1}{n^2} \sum_{i,j=2^{r_n}}^{n} \mathbb{E}[f(x_{\Pi(i)})f(X_{\Pi(j)})].$$

Let us compute the latter expectation, depending on $i = j$ or $i \neq j$. For all $i \in \{2^{r_n}, \ldots, n\}$, $\Pi(i)$ has the uniform distribution on $\mathbb{G}_{r_n}$ so that when $i = j$, $\mathbb{E}[f(X_{\Pi(i)})f(X_{\Pi(j)})] = \mathbb{E}[f^2(X_{\Pi(i)})] = \nu Q^{r_n} f^2$. Let us now treat the case when $i \neq j$. Then $r_n \geq 1$. Independently on $(X, \Pi)$, draw two independent indices $I_{r_n}$ and $J_{r_n}$ uniformly from $\mathbb{G}_{r_n}$. Then since $i \neq j$, the law of $(\Pi(i), \Pi(j))$ is the conditional law of $(I_{r_n}, J_{r_n})$ given $\{I_{r_n} \neq J_{r_n}\}$ so that

$$\begin{aligned}
&\mathbb{E}[f(X_{\Pi(i)})f(X_{\Pi(j)})] \\
&= \mathbb{E}[f(X_{I_{r_n}})f(X_{J_{r_n}})\mathbf{1}_{\{I_{r_n} \neq J_{r_n}\}}]/\mathbb{P}(I_{r_n} \neq J_{r_n}) \\
&= (1 - 2^{-r_n})^{-1} \mathbb{E}[f(X_{I_{r_n}})f(X_{J_{r_n}})\mathbf{1}_{\{I_{r_n} \neq J_{r_n}\}}] \\
&= (1 - 2^{-r_n})^{-1} \left( \mathbb{E}[f(X_{I_{r_n}})f(X_{J_{r_n}})] - \mathbb{E}[f^2(X_{I_{r_n}})\mathbf{1}_{\{I_{r_n} = J_{r_n}\}}] \right) \\
&= (1 - 2^{-r_n})^{-1} \left( \mathbb{E}[f(X_{I_{r_n}})f(X_{J_{r_n}})] - \mathbb{E}[f^2(X_{I_{r_n}})]\mathbb{P}(I_{r_n} = J_{r_n}) \right) \\
&= (1 - 2^{-r_n})^{-1} \left( \mathbb{E}[f(X_{I_{r_n}})f(X_{J_{r_n}})] - 2^{-r_n} \nu Q^{r_n} f^2 \right) \\
&= (1 - 2^{-r_n})^{-1} \sum_{p=0}^{r_n-1} 2^{-p-1} \nu Q^p P(Q^{r_n-p-1} f \otimes Q^{r_n-p-1} f),
\end{aligned}$$

where we have used $\mathbb{P}(I_{r_n} = J_{r_n}) = 2^{-r_n}$ in the second and fifth equalities, the independence of $(X, I_{r_n})$ and $\mathbf{1}_{\{I_{r_n}=J_{r_n}\}}$ in the fourth one, the fact that $X_{I_{r_n}}$ has the same distribution as $Y_{r_n}$ in the fifth one and (7) with $q = r_n$ in the last one. Eventually, we have proved that $b_n^2 = b_n' + b_n''$ with

$$b_n' = \frac{n - 2^{r_n} + 1}{n^2} \nu Q^{r_n} f^2,$$

$$b_n'' = \frac{(n-2^{r_n})(n-2^{r_n}+1)}{n^2(1-2^{-r_n})} \sum_{p=0}^{r_n-1} 2^{-p-1} \nu Q^p P(Q^{r_n-p-1}f \otimes Q^{r_n-p-1}f).$$

Since $(n-2^{r_n}+1)/n^2 \leq 1/n$, and using (ii), (v) and (vi), $\lim_{n\to\infty} b_n' = 0$. As for $b_n''$, let us fix $\varepsilon > 0$ and choose $p_\varepsilon \in \mathbb{N}^*$ such that $2^{-p_\varepsilon} \leq \varepsilon$. From (8), there is a $c_f \geq 0$ such that $|\nu Q^p P(Q^{r_n-p-1}f \otimes Q^{r_n-p-1}f)| \leq c_f$ for all $p$ and $n$ such that $r_n \geq p$. Since $(n-2^{r_n})(n-2^{r_n}+1)/n^2(1-2^{-r_n}) \leq 1$, we then have

$$b_n'' \leq \varepsilon c_f + \sum_{p=0}^{p_\varepsilon-1} |\nu Q^p P(Q^{r_n-p-1}f \otimes Q^{r_n-p-1}f)|.$$

Now, using (10), we get that each term of the latter finite sum tends to 0 as $n$ tends to infinity, so that finally $\lim_{n\to\infty} b_n'' = 0$, which completes the proof. $\square$

PROOF OF THEOREM 12. Considering the function $f - (\mu, Pf)$ leaves us with the case when $(\mu, Pf) = 0$. Let us treat the case of $\overline{M}_{\mathbb{G}_q}(f)$. Observe that $f(\Delta_i) \in L^2$ for all $i \in \mathbb{G}_q$. Indeed, $(z_1, \ldots, z_q)$ denoting the path from the root 1 to $i$ in the tree, $\mathbb{E}[f(\Delta_i)^2] = \nu P_{z_1} \cdots P_{z_q} Pf^2$, which is finite from (iii$'$) and (vi), since $Pf^2 \in F$. Thus, by conditioning on $\mathcal{F}_q$, $\mathbb{E}[f(\Delta_i)f(\Delta_j)] = \mathbb{E}[Pf(X_i)Pf(X_j)]$ for all $i \neq j \in \mathbb{G}_q$, and $\mathbb{E}[f^2(\Delta_i)] = \mathbb{E}[Pf^2(X_i)]$ for all $i \in \mathbb{G}_q$. Hence,

$$\begin{aligned}\mathbb{E}[M_{\mathbb{G}_q}(f)^2] &= \sum_{i,j \in \mathbb{G}_q} \mathbb{E}[f(\Delta_i)f(\Delta_j)] \\ &= \sum_{i \in \mathbb{G}_q} \mathbb{E}[Pf^2(X_i)] + \sum_{i \neq j \in \mathbb{G}_q} \mathbb{E}[Pf(X_i)Pf(X_j)] \\ &= \mathbb{E}[M_{\mathbb{G}_q}(Pf)^2] + \mathbb{E}[M_{\mathbb{G}_q}(Pf^2 - (Pf)^2)],\end{aligned}$$

so that

$$\mathbb{E}[\overline{M}_{\mathbb{G}_q}(f)^2] = \mathbb{E}[\overline{M}_{\mathbb{G}_q}(Pf)^2] + \frac{\mathbb{E}[\overline{M}_{\mathbb{G}_q}(Pf^2 - (Pf)^2)]}{|\mathbb{G}_q|}.$$

We can apply Theorem 11 twice: $\lim_{q\to\infty} \mathbb{E}[\overline{M}_{\mathbb{G}_q}(Pf^2 - (Pf)^2)] = (\mu, Pf^2 - (Pf)^2)$ and $\lim_{q\to\infty} \mathbb{E}[\overline{M}_{\mathbb{G}_q}(Pf)^2] = 0$, so that $\lim_{q\to\infty} \mathbb{E}[\overline{M}_{\mathbb{G}_q}(f)^2] = 0$, that is, $\overline{M}_{\mathbb{G}_q}(f)$ converges to 0 in $L^2$. Using Lemma 13, we extend this result to $\overline{M}_{\mathbb{T}_r}(f)$. The proof for $\overline{M}_n^\Pi(f)$ is similar to Step 3 of the proof of Theorem 11, with the same extra conditioning argument as above. $\square$

2.3. *Strong law of large numbers.* So far, we have proved the weak LLN, that is, convergence in quadratic mean for empirical averages. We now seek for strong LLN. Theorem 14 gives sufficient conditions under which the empirical averages over the $q$th generation and over the first $r+1$ generations converge to a constant with probability one.

THEOREM 14. *Let $F$ satisfy conditions* (i)–(vi). *Let $f \in F$ such that $(\mu, f) = 0$. Assume that there exists $h \in F$ such that*

$$P\left(\sum_{r \in \mathbb{N}} |Q^r f \otimes Q^r f|\right) \leq h.$$

*Then $\overline{M}_{\mathbb{G}_q}(f)$ and $\overline{M}_{\mathbb{T}_r}(f)$ almost surely converge to $0$ as $q \to \infty$.*

PROOF. *Step* 1. Let us first treat $\overline{M}_{\mathbb{G}_q}(f)$. Let us write $\eta_q = \mathbb{E}[\overline{M}_{\mathbb{G}_q}(f)^2]$. It is enough to check that $\sum_{q \in \mathbb{N}} \eta_q < \infty$. Now, using (7) and Fubini's theorem, we have

$$\begin{aligned}
\sum_{q \in \mathbb{N}} \eta_q &= \sum_{q \in \mathbb{N}} \sum_{p=0}^{q} 2^{-p-\mathbf{1}_{\{p<q\}}} \nu Q^p P(Q^{q-p-1} f \otimes Q^{q-p-1} f) \\
&\leq \sum_{p \in \mathbb{N}} \sum_{q=p}^{+\infty} 2^{-p-\mathbf{1}_{\{p<q\}}} |\nu Q^p P(Q^{q-p-1} f \otimes Q^{q-p-1} f)| \\
&\leq \sum_{p \in \mathbb{N}} 2^{-p} \nu Q^p \left(f^2 + P\left(\sum_{r \in \mathbb{N}} |Q^r f \otimes Q^r f|\right)\right) \\
&\leq \sum_{p \in \mathbb{N}} 2^{-p} \nu Q^p (f^2 + h),
\end{aligned}$$

which, from (v) and (vi), is finite, since $f^2 + h \in F$.

*Step* 2. Let us now deal with $\overline{M}_{\mathbb{T}_r}(f)$. From (3), we have $|\overline{M}_{\mathbb{T}_r}(f)| \leq \sum_{q=0}^{r} \frac{|\mathbb{G}_q|}{|\mathbb{T}_r|} |\overline{M}_{\mathbb{G}_q}(f)|$. From Step 1, a.s. $\lim_{q \to \infty} |\overline{M}_{\mathbb{G}_q}(f)| = 0$. It is enough to apply Lemma 13 to get that $\overline{M}_{\mathbb{T}_r}(f)$ a.s. converges to $0$.  □

In particular, we have the following:

COROLLARY 15. *Let $F$ satisfy conditions* (i)–(vi). *Let $f \in F$ such that $(\mu, f) = 0$. Assume there exists $c \in F$ and a nonnegative sequence $(\kappa_r, r \in \mathbb{N})$ such that $\sum_{r \in \mathbb{N}} \kappa_r < \infty$ and*

$$\forall x \in S, \forall r \in \mathbb{N} \qquad |Q^r f(x)| \leq c(x) \kappa_r.$$

*Then $\overline{M}_{\mathbb{G}_q}(f)$ and $\overline{M}_{\mathbb{T}_r}(f)$ almost surely converge to $0$.*

PROOF. From (iii), $h = (\sum_{r\in\mathbb{N}} \kappa_r)P(c \otimes c) \in F$ and $P(\sum_{r\in\mathbb{N}} |Q^r f \otimes Q^r f|) \leq h$, so that Theorem 14 gives the result. □

REMARK 16. In the case when $\kappa_r = \kappa^r$ for some $\kappa \in (0,1)$, we speak of geometric ergodicity. Geometric ergodicity implies the almost sure convergence of $\overline{M}_{\mathbb{G}_q}(f)$ and $\overline{M}_{\mathbb{T}_r}(f)$.

REMARK 17. Assume that the state space $S$ is finite, and that the induced MC $Y$ is irreducible and aperiodic. Then $Y$ is ergodic and, $\mu$ standing for its stationary distribution, the sequence of functions $(Q^r g, r \in \mathbb{N})$ uniformly converges to $(\mu, g)$ with exponential speed. Taking $F = \mathcal{B}(S)$ and $f = g - (\mu, g)$, Corollary 15 applies: $\overline{M}_{\mathbb{G}_q}(f)$ and $\overline{M}_{\mathbb{T}_r}(f)$ almost surely converge to 0, that is, $\overline{M}_{\mathbb{G}_q}(g)$ and $\overline{M}_{\mathbb{T}_r}(g)$ almost surely converge to $(\mu, g)$. This covers the main result in [18, 23, 24] when applied to the binary tree.

In the case when $f$ depends on the mother–daughters triangle $\Delta_n = (X_n, X_{2n}, X_{2n+1})$, we can prove as well the following:

THEOREM 18. Let $F$ satisfy conditions (i)–(vi). Let $f \in \mathcal{B}(S^3)$ such that $Pf$ and $Pf^2$ exist and belong to $F$, with $(\mu, Pf) = 0$.

(i) Assume that there exists $h \in F$ such that $P(\sum_{r\in\mathbb{N}} |Q^r Pf \otimes Q^r Pf|) \leq h$. Then $\overline{M}_{\mathbb{G}_q}(f)$ and $\overline{M}_{\mathbb{T}_r}(f)$ almost surely converge to 0.

(ii) In particular, if there exists $c \in F$ and a nonnegative sequence $(\kappa_r, r \in \mathbb{N})$ such that $\sum_{r\in\mathbb{N}} \kappa_r < \infty$ and

$$\forall x \in S, \forall r \in \mathbb{N} \quad |Q^r Pf(x)| \leq c(x)\kappa_r,$$

then $\overline{M}_{\mathbb{G}_q}(f)$ and $\overline{M}_{\mathbb{T}_r}(f)$ almost surely converge to 0.

2.4. *Central limit theorem.* We are now interested in proving a central limit theorem (CLT) for the $\mathbb{T}$-MC $(X_n)$. This will be done by using a CLT for martingales.

THEOREM 19. Let $F$ satisfy (i)–(vi). Let $f \in \mathcal{B}(S^3)$ such that $Pf^2$ and $Pf^4$ exist and belong to $F$. Assume that $Pf = 0$. Then $n^{-1/2}M_n^\Pi(f)$ converges in distribution to the Gaussian law $\mathcal{N}(0, s^2)$, where $s^2 = (\mu, Pf^2)$.

PROOF. Let $M_0^\Pi(f) = 0$, $\mathcal{H}_0 = \sigma(X_1)$ and $\mathcal{H}_n = \sigma(\Delta_{\Pi(i)}, \Pi(i+1), 1 \leq i \leq n)$ for all $n \geq 1$. Note that $X_{\Pi(i)}$ is $\mathcal{H}_{i-1}$-measurable and that, conditionally on $\mathcal{H}_{i-1}$, $\Delta_{\Pi(i)}$ has distribution $\delta_{X_{\Pi(i)}} \otimes P(X_{\Pi(i)}, \cdot)$. Since $Pf = 0$, $(M_n^\Pi(f), n \geq 0)$ is a $(\mathcal{H}_n)$-martingale. It has bracket

$$\langle M^\Pi(f)\rangle_n = \sum_{i=1}^n \mathbb{E}[f^2(\Delta_{\Pi(i)})|\mathcal{H}_{i-1}] = \sum_{i=1}^n Pf^2(X_{\Pi(i)}) = M_n^\Pi(Pf^2).$$

According to Theorem 11, since $Pf^2 \in F$, the sequence $n^{-1}\langle M^\Pi(f)\rangle_n$ converges to $(\mu, Pf^2) = s^2$ in $L^2$, and thus, in probability. It remains to check Liapunov's condition, say, for the fourth moment, that is, to prove that the sequence of positive r.v. $(L_n, n \geq 1)$ defined by

$$L_n = \frac{1}{n^2} \sum_{i=1}^{n} \mathbb{E}[f^4(\Delta_{\Pi(i)})|\mathcal{H}_{i-1}]$$

tends in probability to 0 (see, e.g., [10]). But $L_n = \overline{M}_n^\Pi(Pf^4)/n$ and $\overline{M}_n^\Pi(Pf^4)$ converges to $(\mu, Pf^4)$ in quadratic mean, so that $L_n$ converges to 0 in probability. □

In the general case when $Pf \neq 0$, we have:

COROLLARY 20. *Let $F$ satisfy* (i)–(vi). *Let $f \in \mathcal{B}(S^3)$ such that $Pf$, $Pf^2$ and $Pf^4$ exist and belong to $F$. Then $n^{-1/2}(M_n^\Pi(f) - M_n^\Pi(Pf))$ converges in distribution to $\mathcal{N}(0, s^2)$, where $s^2 = (\mu, Pf^2) - (\mu, (Pf)^2)$.*

PROOF. It is enough to apply Theorem 19 to the function $g$ defined by $g(x, y, z) = f(x, y, z) - Pf(x)$. □

Considering the subsequence of indices $n = |\mathbb{T}_r|$, $r \in \mathbb{N}$, we can state the following:

COROLLARY 21. *Let $F$ satisfy* (i)–(vi). *Let $f \in \mathcal{B}(S^3)$ such that $Pf$, $Pf^2$ and $Pf^4$ exist and belong to $F$. Then $|\mathbb{T}_r|^{1/2}(\overline{M}_{\mathbb{T}_r}(f) - \overline{M}_{\mathbb{T}_r}(Pf))$ converges in distribution to $\mathcal{N}(0, s^2)$, where $s^2 = (\mu, Pf^2) - (\mu, (Pf)^2)$.*

If we take $F$ to be $\mathcal{C}_b(S)$, we get the following:

COROLLARY 22. *Assume that the induced MC $(Y_r, r \in \mathbb{N})$ is ergodic, with stationary distribution $\mu$ (see Definition 7). Then, for any $f \in \mathcal{C}_b(S^3)$, $n^{-1/2}(M_n^\Pi(f) - M_n^\Pi(Pf))$ converges in distribution to $\mathcal{N}(0, s^2)$, where $s^2 = (\mu, Pf^2) - (\mu, (Pf)^2)$.*

REMARK 23. Note that the normalizing factor is the square root of the number of individuals, $n$ or $|\mathbb{T}_r|$, and not the square root of the number of generations, $r_n$ or $r$, as one might have thought. Convergence is fast with $r$: with 20 generations ($r = 19$), the normalizing factor $|\mathbb{T}_r|^{1/2}$ is approximately $10^3$.

Using characteristic functions, it is easy to generalize Corollary 22 to the case when $f$ is vector-valued:

COROLLARY 24. *Let $F$ satisfy* (i)–(vi). *Let $f_1, \ldots, f_d \in \mathcal{B}(S^3)$ such that $Pf_i$, $P(f_if_j)$ and $P(f_if_jf_kf_l)$ exist and belong to $F$ for all $i,j,k,l$. Let $f = (f_1, \ldots, f_d)$. Then $n^{-1/2}(M_n^\Pi(f) - M_n^\Pi(Pf))$ converges in distribution to the $d$-dimensional Gaussian law $\mathcal{N}_d(0, \Sigma)$, where $\Sigma_{ij} = (\mu, P(f_if_j)) - (\mu, Pf_iPf_j)$.*

## 3. Detection of cellular aging.

3.1. *Limit theorems in the BAR model* (1). Here we seek to apply the results in Section 2 to model (1).

3.1.1. *Weak law of large numbers and central limit theorem.* In this section we take $F$ to be the set $C_{\mathrm{pol}}(\mathbb{R})$ of continuous and polynomially growing functions, that is, the set of all continuous functions $f : \mathbb{R} \to \mathbb{R}$ such that there exists $c \geq 0$ and $m \in \mathbb{N}$ such that, for all $x \in \mathbb{R}$,

$$|f(x)| \leq c(1 + |x|^m).$$

In order to apply Theorems 11 and 12 and Corollary 24, we need to check conditions (i)–(vi). Conditions (i) and (ii) are obvious. The next lemma states that condition (iii) is fulfilled too.

LEMMA 25. *Let $f, g \in C_{\mathrm{pol}}(\mathbb{R})$. Then $f \otimes g \in L^1(P(x, \cdot))$ and $P(f \otimes g) \in C_{\mathrm{pol}}(\mathbb{R})$.*

PROOF. Let $G_0, G_1$ be two independent standard Gaussian variables. Let

$$G = \begin{pmatrix} G_0 \\ G_1 \end{pmatrix}, \qquad M(x) = \begin{pmatrix} \alpha_0 x + \beta_0 \\ \alpha_1 x + \beta_1 \end{pmatrix} \quad \text{and} \quad \Lambda = \sigma \begin{pmatrix} 1 & 0 \\ \rho & \sqrt{1 - \rho^2} \end{pmatrix}.$$

Then $G(x) \equiv M(x) + \Lambda G$ has distribution $P(x, \cdot)$. Hence, $P(|f \otimes g|)(x) = \mathbb{E}[|f \otimes g(G(x))|]$. Now, $\|\ \|$ denoting the Euclidian norm in $\mathbb{R}^2$, we can find $c \geq 0$ and $m \in \mathbb{N}$ such that $|f \otimes g(G(x))| \leq c(1 + \|G(x)\|^m)$. Since $\|G(x)\| \leq c(1 + |x| + \|G\|)$ for a constant $c$ and $\mathbb{E}[\|G\|^m] < \infty$, we may eventually find a $c \geq 0$ such that, for all $x \in \mathbb{R}$,

$$P(|f \otimes g|)(x) \leq c(1 + |x|^m),$$

which completes the proof. □

LEMMA 26. *$C_{\mathrm{pol}}(\mathbb{R})$ fulfills conditions* (iv)–(v); *$\mu$ is the stationary distribution of $Y$.*

PROOF. In the BAR model (1), the induced MC has the stochastic dynamics

(11) $$Y_{r+1} = a_{r+1} Y_r + b_{r+1},$$

where $((a_r, b_r), r \in \mathbb{N}^*)$ is a sequence of i.i.d. r.v., independent of $Y_0$. Precisely, $a_{r+1} = \alpha_{\zeta_{r+1}}$, $b_{r+1} = \beta_{\zeta_{r+1}} + \varepsilon'_{r+1}$, where $(\varepsilon'_q, q \in \mathbb{N}^*)$ and $(\zeta_q, q \in \mathbb{N}^*)$ are independent sequences of i.i.d. r.v., independent of $Y_0$, each $\varepsilon'_q$ has law $\mathcal{N}(0, \sigma^2)$ and each $\zeta_q$ is a balanced Bernoulli r.v., that is, $\mathbb{P}(\zeta_q = 0) = \mathbb{P}(\zeta_q = 1) = 1/2$. Bougerol and Picard [7] call the sequence $(Y_r, r \in \mathbb{N})$ a generalized autoregressive process. It is often simply called AR(1), for AutoRegressive of order 1, in the literature. We have

$$Y_r = a_r a_{r-1} \cdots a_2 a_1 Y_0 + \sum_{k=1}^{r} a_r a_{r-1} \cdots a_{k+1} b_k.$$

Since the r.v. $((a_r, b_r), r \in \mathbb{N}^*)$ are i.i.d., $Y_r$ has the same distribution as

(12) $$Z_r = a_1 a_2 \cdots a_{r-1} a_r Y_0 + \sum_{k=1}^{r} a_1 a_2 \cdots a_{k-1} b_k.$$

Let us first prove (v). Let $f \in C_{\text{pol}}(\mathbb{R})$ and $x \in \mathbb{R}$, and let us denote $S = \sum_{k=1}^{\infty} |a_1 a_2 \cdots a_{k-1} b_k|$. From (12), $|Z_r| \leq |Y_0| + S$ for all $r \in \mathbb{N}$, so that we can find $c \geq 0$ and $m \in \mathbb{N}^*$ such that $\mathbb{E}_x[|f(Z_r)|] \leq c(|x|^m + \mathbb{E}[S^m])$. Now, let us denote $\alpha = \max\{|\alpha_0|, |\alpha_1|\} < 1$. Using the triangle inequality in the first line and the fact that the $L^m$-norm of $b_k$, $\|b_k\|_{L^m} \equiv c_m$, does not depend on $k$ in the last one,

(13) $$\|S\|_{L^m} \leq \sum_{k=1}^{\infty} \|a_1 a_2 \cdots a_{k-1} b_k\|_{L^m}$$
$$\leq \sum_{k=1}^{\infty} \alpha^{k-1} \|b_k\|_{L^m} = c_m \sum_{k=1}^{\infty} \alpha^{k-1} < \infty.$$

Eventually, $|Q^r f(x)| \leq \mathbb{E}_x[|f(Y_r)|] = \mathbb{E}_x[|f(Z_r)|] \leq c'_m (1 + |x|^m)$ for some $c'_m$ which does not depend on $r$, which proves (v).

Let us now prove (iv). Since $|a_1 a_2 \cdots a_{r-1} a_r Y_0| \leq \alpha^r |Y_0|$, we have that a.s. $\lim_{r \to \infty} a_1 a_2 \cdots a_{r-1} a_r Y_0 = 0$. Besides, the sum in (12) a.s. converges as $r$ grows to infinity as $\mathbb{E}[|S|^m] < \infty$. Eventually, the sequence $(Z_r, r \in \mathbb{N})$ almost surely converges to

(14) $$Z_\infty = \sum_{k=1}^{\infty} a_1 a_2 \cdots a_{k-1} b_k.$$

Let $\mu$ denote the distribution of $Z_\infty$. Then $C_{\text{pol}}(\mathbb{R}) \subset L^1(\mu)$. Indeed, $\|Z_\infty\|_{L^m} \leq \|S\|_{L^m} < \infty$ for all $m \in \mathbb{N}$. Let us eventually prove that $\lim_{r \to \infty} \mathbb{E}_x[f(Y_r)] =$

$(\mu, f)$. Since $|f(Z_r)| \le c(|Y_0|^m + S^m) \in L^1(\mathbb{P}_x)$ and $(Z_r, r \in \mathbb{N})$ almost surely converges to $Z_\infty$, we can apply Lebesgue's dominated convergence theorem and get that

$$\lim_{r \to \infty} \mathbb{E}_x[f(Y_r)] = \lim_{r \to \infty} \mathbb{E}_x[f(Z_r)] = \mathbb{E}_x[f(Z_\infty)] = (\mu, f).$$

Condition (iv) is now fully checked, and $\mu$ is the unique stationary distribution of $Y$. □

Let us denote by $C_{\text{pol}}(\mathbb{R}^3)$ the set of all continuous and polynomially growing functions $f : \mathbb{R}^3 \to \mathbb{R}$. Since $C_{\text{pol}}(\mathbb{R}^3)^2 \subset C_{\text{pol}}(\mathbb{R}^3)$ and $P(C_{\text{pol}}(\mathbb{R}^3)) \subset C_{\text{pol}}(\mathbb{R})$, Theorems 11 and 12 and Corollary 24 imply:

PROPOSITION 27. *In the BAR model* (1), *assume that the distribution of the ancestor $X_1$, $\nu$, has finite moments of all orders. Let $\mu$ be the unique stationary distribution of the induced MC $(Y_r, r \in \mathbb{N})$. Then:*

(i) *for all $f \in C_{\text{pol}}(\mathbb{R})$, $\overline{M}_{\mathbb{G}_q}(f)$, $\overline{M}_{\mathbb{T}_r}(f)$ and $\overline{M}_n^\Pi(f)$ converge to $(\mu, f)$ in $L^2$,*

(ii) *for all $f \in C_{\text{pol}}(\mathbb{R}^3)$, $\overline{M}_{\mathbb{G}_q}(f)$, $\overline{M}_{\mathbb{T}_r}(f)$ and $\overline{M}_n^\Pi(f)$ converge to $(\mu, Pf)$ in $L^2$,*

(iii) *for all $f_1, \ldots, f_d \in C_{\text{pol}}(\mathbb{R}^3)$, $n^{-1/2}(M_n^\Pi(f) - M_n^\Pi(Pf))$ converges in distribution to $\mathcal{N}_d(0, \Sigma)$, where $f = (f_1, \ldots, f_d)$ and $\Sigma_{ij} = (\mu, P(f_i f_j)) - (\mu, Pf_i Pf_j)$.*

3.1.2. *Strong law of large numbers.* We can also derive almost sure convergence results:

PROPOSITION 28. *With the assumptions of Proposition 27:*

(i) *almost surely, for any $f \in C_{\text{pol}}(\mathbb{R})$, $\overline{M}_{\mathbb{G}_q}(f)$ and $\overline{M}_{\mathbb{T}_r}(f)$ converge to $(\mu, f)$,*

(ii) *almost surely, for any $f \in C_{\text{pol}}(\mathbb{R}^3)$, $\overline{M}_{\mathbb{G}_q}(f)$ and $\overline{M}_{\mathbb{T}_r}(f)$ converge to $(\mu, Pf)$.*

PROOF. Let us take $F$ to be the set $C^1_{\text{pol}}(\mathbb{R})$ of all $C^1$ functions $f : \mathbb{R} \to \mathbb{R}$ such that $|f| + |f'|$ is bounded above by a polynomial. One can easily check that $C^1_{\text{pol}}(\mathbb{R})$ satisfies (i)–(v).

*Step* 1. Let us first prove that,

$$(15) \quad \forall f \in C^1_{\text{pol}}(\mathbb{R}) \qquad \mathbb{P}\left(\lim_{q \to \infty} \overline{M}_{\mathbb{G}_q}(f) = \lim_{r \to \infty} \overline{M}_{\mathbb{T}_r}(f) = (\mu, f)\right) = 1.$$

Let $f \in C^1_{\text{pol}}(\mathbb{R})$. We want to apply Corollary 15 with $F = C^1_{\text{pol}}(\mathbb{R})$ and to the function $g = f - (\mu, f) \in C^1_{\text{pol}}(\mathbb{R})$. First note that $Q^r g(x) = Q^r f(x) - (\mu, f) =$

$\mathbb{E}_x[f(Z_r) - f(Z_\infty)]$, so that using Cauchy–Schwarz's inequality,

$$|Q^r g(x)| \leq \mathbb{E}_x[W_r|Z_r - Z_\infty|] \leq (\mathbb{E}_x[W_r^2]\mathbb{E}_x[(Z_r - Z_\infty)^2])^{1/2},$$

where $W_r = \sup_{z \in [Z_r, Z_\infty]} |f'(z)|$. We can find $c_1 \geq 0$ and $m \in \mathbb{N}$ such that for all $z \in \mathbb{R}$, $|f'(z)|^2 \leq c_1(1 + |z|^m)$ so that, using (12) and (13), there is a $c'_m \geq 0$ such that, for all $x \in \mathbb{R}$ and $r \in \mathbb{N}$,

$$\mathbb{E}_x[W_r^2] \leq \mathbb{E}_x\left[\sup_{z \in [Z_r, Z_\infty]} c_1(1 + |z|^m)\right]$$

$$\leq c_1(1 + \mathbb{E}_x[|Z_r|^m] + \mathbb{E}_x[|Z_\infty|^m]) \leq c'_m(1 + |x|^m).$$

Moreover, $Z_r - Z_\infty = a_1 a_2 \cdots a_{r-1} a_r Y_0 - \sum_{k=r+1}^{+\infty} a_1 a_2 \cdots a_{k-1} b_k$ so that $\|Z_r - Z_\infty\|_{L^2(\mathbb{P}_x)} \leq \alpha^r x + \left\|\sum_{k=r+1}^{+\infty} a_1 a_2 \cdots a_{k-1} b_k\right\|_{L^2}$ ($\alpha = \max\{|\alpha_0|, |\alpha_1|\}$). Now we have $\left\|\sum_{k=r+1}^{+\infty} a_1 a_2 \cdots a_{k-1} b_k\right\|_{L^2} \leq c_2 \sum_{k=r+1}^{+\infty} \alpha^{k-1} = c_2 \alpha^r/(1-\alpha)$, where $c_2 = \|b_k\|_{L^2}$ does not depend on $k$. Thus, we can find $c_3 \geq 0$ such that for all $x \in \mathbb{R}$ and $r \in \mathbb{N}$,

$$\mathbb{E}_x[(Z_r - Z_\infty)^2] \leq c_3 \alpha^{2r}(1 + x^2).$$

Eventually,

(16) $\qquad |Q^r g(x)| \leq (c'_m c_3 \alpha^{2r}(1 + |x|^m)(1 + x^2))^{1/2} \leq c(x)\kappa_r,$

with $\kappa_r = \alpha^r$ and a function $c \in C^1_{\text{pol}}(\mathbb{R})$. Corollary 15 implies that $\overline{M}_{\mathbb{G}_q}(g)$ and $\overline{M}_{\mathbb{T}_r}(g)$ almost surely converge to 0, that is, $\overline{M}_{\mathbb{G}_q}(f)$ and $\overline{M}_{\mathbb{T}_r}(f)$ almost surely converge to $(\mu, f)$, which proves (15).

*Step* 2. Let us now prove that, almost surely, the empirical distributions $\overline{M}_{\mathbb{G}_q}$ and $\overline{M}_{\mathbb{T}_r}$ weakly converge to $\mu$. There exists a sequence $(f_p, p \in \mathbb{N})$ of $C^\infty$ functions with compact support which characterizes convergence in distribution. Hence, it is enough to show that, almost surely, for all $p \in \mathbb{N}$, $\lim_{q \to \infty} \overline{M}_{\mathbb{G}_q}(f_p) = \lim_{r \to \infty} \overline{M}_{\mathbb{T}_r}(f_p) = (\mu, f_p)$. But this immediately follows from Step 1, since $f_p \in C^1_{\text{pol}}(\mathbb{R})$.

*Step* 3. Let us now prove assertion (i). Let us deal with $\overline{M}_{\mathbb{G}_q}$ (the proof for $\overline{M}_{\mathbb{T}_r}$ is similar). For $k, l \in \mathbb{N}$, let us write $f_{k,l}(x) = k(1 + x^{2l})$. Since $f_{k,l} \in C^1_{\text{pol}}(\mathbb{R})$, from Step 1, almost surely,

(17) $\qquad \forall k, l \in \mathbb{N} \quad \lim_{q \to \infty} \overline{M}_{\mathbb{G}_q}(f_{k,l}) = (\mu, f_{k,l}).$

From Step 2, the empirical distributions $\overline{M}_{\mathbb{G}_q}$ a.s. weakly converge to $\mu$. Besides, for all $f \in C_{\text{pol}}(\mathbb{R})$, there exists $k, l \in \mathbb{N}$ such that $f^2 \leq f_{k,l}$. Thus, from (17), a.s. for all $f \in C_{\text{pol}}(\mathbb{R})$, the sequence $(\overline{M}_{\mathbb{G}_q}(f^2), q \in \mathbb{N})$ is bounded, which proves that a.s. every $f \in C_{\text{pol}}(\mathbb{R})$ is $(\overline{M}_{\mathbb{G}_q}, q \in \mathbb{N})$-uniformly integrable. Hence, a.s., for all $f \in C_{\text{pol}}(\mathbb{R})$, $\overline{M}_{\mathbb{G}_q}(f)$ converges to $(\mu, f)$.

*Step* 4. The proof of (ii) is similar to the proof of (i). □

REMARK 29. A natural choice for $\nu$ is the stationary distribution $\mu$. Indeed, the ancestor $X_1$ is picked from a metacolony that has evolved for a long time, so that in the BAR model (1) its distribution should be close to $\mu$. With this particular choice, we can apply Propositions 27 and 28. Indeed, $\mu$ has finite moments of all orders, since $C_{\text{pol}}(\mathbb{R}) \subset L^1(\mu)$.

3.2. *Estimation of the parameters.* We seek to estimate the 4-dimensional parameter $\theta = (\alpha_0, \beta_0, \alpha_1, \beta_1)$, as well as $\sigma^2$ and $\rho$. Assume we observe a complete subtree $\mathbb{T}_{r+1}$. Then, since the couples $(\varepsilon_{2i}, \varepsilon_{2i+1})$ are i.i.d. bivariate Gaussian vectors, the maximum likelihood estimator $\hat{\theta}^r = (\hat{\alpha}_0^r, \hat{\beta}_0^r, \hat{\alpha}_1^r, \hat{\beta}_1^r)$ of $\theta$ is also the least squares one: for $\varepsilon \in \{0, 1\}$,

$$\begin{cases} \hat{\alpha}_\varepsilon^r = \dfrac{|\mathbb{T}_r|^{-1} \sum_{i \in \mathbb{T}_r} X_i X_{2i+\varepsilon} - (|\mathbb{T}_r|^{-1} \sum_{i \in \mathbb{T}_r} X_i)(|\mathbb{T}_r|^{-1} \sum_{i \in \mathbb{T}_r} X_{2i+\varepsilon})}{|\mathbb{T}_r|^{-1} \sum_{i \in \mathbb{T}_r} X_i^2 - (|\mathbb{T}_r|^{-1} \sum_{i \in \mathbb{T}_r} X_i)^2}, \\ \hat{\beta}_\varepsilon^r = |\mathbb{T}_r|^{-1} \sum_{i \in \mathbb{T}_r} X_{2i+\varepsilon} - \hat{\alpha}_0^r |\mathbb{T}_r|^{-1} \sum_{i \in \mathbb{T}_r} X_i. \end{cases}$$

Hence, $\hat{\alpha}_0^r$ (resp. $\hat{\alpha}_1^r$) is the empirical correlation between new (resp. old) pole daughters and their mothers. We shall denote by $\mathbf{xy}$ (resp. $\mathbf{x}, \mathbf{y}, \mathbf{x}^2$) the element of $C_{\text{pol}}(\mathbb{R}^3)$ defined by $(x, y, z) \mapsto xy$ (resp. $x, y, x^2$).

REMARK 30. Note that $(\mu, \mathbf{x}^2) - (\mu, \mathbf{x})^2 > 0$. Indeed, it is nonnegative, and if it were 0, $\mu$ would be a Dirac mass. Now a Dirac mass cannot be stationary for $Y$ because $\sigma^2 > 0$.

PROPOSITION 31. *In the BAR model* (1), *assume that the distribution of the ancestor $X_1$, $\nu$, has finite moments of all orders. Then $(\hat{\theta}^r, r \in \mathbb{N})$ is a strongly consistent estimator of $\theta$.*

PROOF. Let us treat $\hat{\alpha}_0^r$. Convergence of $\hat{\beta}_0^r$, $\hat{\alpha}_1^r$ and $\hat{\beta}_1^r$ may be treated in a similar way. Note that $\hat{\alpha}_0^r = C_r/B_r$ with

$$C_r = \overline{M}_{\mathbb{T}_r}(\mathbf{xy}) - \overline{M}_{\mathbb{T}_r}(\mathbf{x})\overline{M}_{\mathbb{T}_r}(\mathbf{y}) \quad \text{and} \quad B_r = \overline{M}_{\mathbb{T}_r}(\mathbf{x}^2) - \overline{M}_{\mathbb{T}_r}(\mathbf{x})^2.$$

Since $P(\mathbf{xy})(x) = x(\alpha_0 x + \beta_0)$ and $P(\mathbf{y})(x) = \alpha_0 x + \beta_0$, Proposition 28 implies that $C_r$ a.s. converges to $(\mu, \mathbf{x}(\alpha_0 \mathbf{x} + \beta_0)) - (\mu, \mathbf{x})(\mu, \alpha_0 \mathbf{x} + \beta_0) = \alpha_0((\mu, \mathbf{x}^2) - (\mu, \mathbf{x})^2)$ and $B_r$ a.s. converges to $(\mu, \mathbf{x}^2) - (\mu, \mathbf{x})^2$, which from Remark 30 is positive, so that $\hat{\alpha}_0^r$ a.s. converges to $\alpha_0$. □

REMARK 32. Let us denote $\bar{\alpha} = (\alpha_0 + \alpha_1)/2$, $\bar{\beta} = (\beta_0 + \beta_1)/2$ and so on. Then

$$(18) \qquad (\mu, \mathbf{x}) = \frac{\bar{\beta}}{1 - \bar{\alpha}} \quad \text{and} \quad (\mu, \mathbf{x}^2) = \frac{2\overline{\alpha\beta}\bar{\beta}/(1 - \bar{\alpha}) + \overline{\beta^2} + \sigma^2}{1 - \overline{\alpha^2}}.$$

Indeed, recalling (11) and (14), $Z_\infty$ has the same law as $a_1 Z_\infty + b_1$, where the pair $(a_1, b_1)$ is independent of $Z_\infty$ and takes values $(\alpha_0, \beta_0)$ and $(\alpha_1, \beta_1)$ with probability 1/2. Hence, $(\mu, \mathbf{x}) = \mathbb{E}[Z_\infty] = \mathbb{E}[a_1 Z_\infty + b_1] = \bar{\alpha}(\mu, \mathbf{x}) + \bar{\beta}$, as announced. Likewise,

$$(\mu, \mathbf{x}^2) = \mathbb{E}[Z_\infty^2] = \mathbb{E}[(a_1 Z_\infty + b_1)^2] = \mathbb{E}[a_1^2 Z_\infty^2] + 2\mathbb{E}[a_1 b_1 Z_\infty] + \mathbb{E}[b_1^2]$$
$$= \overline{\alpha^2}(\mu, \mathbf{x}^2) + 2\overline{\alpha\beta}(\mu, \mathbf{x}) + \overline{\beta^2} + \sigma^2,$$

from which we deduce the second equality in (18).

From the preceding remark, we define two continuous functions $\mu_1 : \Theta \to \mathbb{R}$ and $\mu_2 : \Theta \times \mathbb{R}_+^* \to \mathbb{R}$ by writing

(19) $$(\mu, \mathbf{x}) = \mu_1(\theta) \quad \text{and} \quad (\mu, \mathbf{x}^2) = \mu_2(\theta, \sigma^2),$$

where $\theta = (\alpha_0, \beta_0, \alpha_1, \beta_1) \in \Theta = (-1, 1) \times \mathbb{R} \times (-1, 1) \times \mathbb{R}$. Let us now build a confidence region for $\theta$.

PROPOSITION 33. *In the BAR model* (1), *assume that the distribution $\nu$ of the ancestor $X_1$ has finite moments of all orders. Let $\mu$ be the unique stationary distribution of the induced MC $(Y_r, r \in \mathbb{N})$. Then $|\mathbb{T}_r|^{1/2}(\hat{\theta}^r - \theta)$ converges in law to $\mathcal{N}_4(0, \Sigma')$, where*

(20)
$$\Sigma' = \sigma^2 \begin{pmatrix} K & \rho K \\ \rho K & K \end{pmatrix}$$

$$\text{with } K = \frac{1}{\mu_2(\theta, \sigma^2) - \mu_1(\theta)^2} \begin{pmatrix} 1 & -\mu_1(\theta) \\ -\mu_1(\theta) & \mu_2(\theta, \sigma^2) \end{pmatrix}.$$

PROOF. For $f_1, \ldots, f_d \in C_{\text{pol}}(\mathbb{R}^3)$, we denote $f = (f_1, \ldots, f_d)$ and $U^r(f) = |\mathbb{T}_r|^{1/2}(M_{\mathbb{T}_r}(f) - M_{\mathbb{T}_r}(Pf))$. Let us denote $\zeta_r = |\mathbb{T}_r|^{1/2}(\hat{\theta}^r - \theta)$. We first observe that $\zeta_r = \varphi(U^r(f), A_r, B_r)$ with $f = (\mathbf{xy}, \mathbf{y}, \mathbf{xz}, \mathbf{z})$, $\varphi(u, a, b) = M(a, b)u$,

$$M(a, b) = \begin{pmatrix} 1/b & -a/b & 0 & 0 \\ -a/b & (b+a^2)/b & 0 & 0 \\ 0 & 0 & 1/b & -a/b \\ 0 & 0 & -a/b & (b+a^2)/b \end{pmatrix},$$

$A_r = \overline{M}_{\mathbb{T}_r}(\mathbf{x})$ and $B_r = \overline{M}_{\mathbb{T}_r}(\mathbf{x}^2) - \overline{M}_{\mathbb{T}_r}(\mathbf{x})^2$. From Proposition 27(iii), $U^r(f)$ converges in distribution to $G \sim \mathcal{N}_4(0, \Sigma)$ with

$$\Sigma = \sigma^2 \begin{pmatrix} \mu_2(\theta, \sigma^2) & \mu_1(\theta) & \rho\mu_2(\theta, \sigma^2) & \rho\mu_1(\theta) \\ \mu_1(\theta) & 1 & \rho\mu_1(\theta) & \rho \\ \rho\mu_2(\theta, \sigma^2) & \rho\mu_1(\theta) & \mu_2(\theta, \sigma^2) & \mu_1(\theta) \\ \rho\mu_1(\theta) & \rho & \mu_1(\theta) & 1 \end{pmatrix}.$$

Besides, Proposition 27(i) implies that $(A_r, B_r)$ converges in law to the constant $(a, b) = (\mu_1(\theta), \mu_2(\theta, \sigma^2) - \mu_1(\theta)^2)$. Thus, Slutsky's theorem states that $(U^r(f), A_r, B_r)$ converges in law to $(G, a, b)$. Then, by continuity of $\varphi$ on $\mathbb{R} \times \mathbb{R} \times \mathbb{R}_+^*$ and recalling from Remark 30 that $b > 0$, $\zeta_r = \varphi(U^r(f), A_r, B_r)$ converges in law to $\varphi(G, a, b) = M(a, b)G$, which is a centered Gaussian vector with covariance matrix $\Sigma' = M(a, b)\Sigma M(a, b)^t$. Now we have

$$\Sigma = \sigma^2 \begin{pmatrix} L & \rho L \\ \rho L & L \end{pmatrix} \quad \text{with } L = \begin{pmatrix} \mu_2(\theta, \sigma^2) & \mu_1(\theta) \\ \mu_1(\theta) & 1 \end{pmatrix},$$

$$M(a, b) = \begin{pmatrix} K & 0 \\ 0 & K \end{pmatrix}.$$

Since $LK = I_2$, the $2 \times 2$ identity matrix, we get

$$\Sigma' = \sigma^2 \begin{pmatrix} K & 0 \\ 0 & K \end{pmatrix} \begin{pmatrix} L & \rho L \\ \rho L & L \end{pmatrix} \begin{pmatrix} K & 0 \\ 0 & K \end{pmatrix} = \sigma^2 \begin{pmatrix} K & \rho K \\ \rho K & K \end{pmatrix}$$

which completes the proof. □

We also need to estimate the conditional variance, $\sigma^2$, and the conditional sister-sister correlation, $\rho$. Since $\sigma^2$ is the common expectation of the i.i.d. r.v. $(\varepsilon_i^2, i \geq 2)$, it is naturally estimated, given a complete observation $(X_i, i \in \mathbb{T}_{r+1})$, by

$$\hat{\sigma}_r^2 = \frac{1}{2|\mathbb{T}_r|} \sum_{i \in \mathbb{T}_r} (\hat{\varepsilon}_{2i}^2 + \hat{\varepsilon}_{2i+1}^2),$$

where

$$\begin{cases} \hat{\varepsilon}_{2n} = X_{2n} - \hat{\alpha}_0^r X_n - \hat{\beta}_0^r, \\ \hat{\varepsilon}_{2n+1} = X_{2n+1} - \hat{\alpha}_1^r X_n - \hat{\beta}_1^r, \end{cases}$$

are the residues. Likewise, since $\rho = \text{Cov}(\varepsilon_{2i}, \varepsilon_{2i+1})/\sigma^2$, it is naturally estimated by

$$\hat{\rho}_r = \frac{1}{\hat{\sigma}_r^2 |\mathbb{T}_r|} \sum_{i \in \mathbb{T}_r} \hat{\varepsilon}_{2i} \hat{\varepsilon}_{2i+1}.$$

We have checked that $(\hat{\sigma}_r^2, \hat{\rho}_r)$ is the maximum likelihood estimator of $(\sigma^2, \rho)$.

PROPOSITION 34. *In the BAR model* (1), *assume that the distribution of the ancestor* $X_1$, $\nu$, *has finite moments of all orders. Then* $((\hat{\sigma}_r^2, \hat{\rho}_r), r \in \mathbb{N})$ *is a strongly consistent estimator of* $(\sigma^2, \rho)$.

PROOF. Let us first deal with $\hat{\sigma}_r^2$. Observe that

$$\hat{\sigma}_r^2 = \frac{1}{2|\mathbb{T}_r|} \sum_{i \in \mathbb{T}_r} f(\Delta_i, \hat{\theta}^r),$$

where $f(\Delta, \theta) = (y - \alpha_0 x - \beta_0)^2 + (z - \alpha_1 x - \beta_1)^2$, with $\Delta = (x, y, z)$ and $\theta = (\alpha_0, \beta_0, \alpha_1, \beta_1)$. Thus, we have $\hat{\sigma}_r^2 = \overline{M}_{\mathbb{T}_r}(f(\cdot, \theta))/2 + D_r$, with

$$D_r = \frac{1}{2|\mathbb{T}_r|} \sum_{i \in \mathbb{T}_r} (f(\Delta_i, \hat{\theta}^r) - f(\Delta_i, \theta)).$$

Since $f(\cdot, \theta) \in C_{\text{pol}}(\mathbb{R}^3)$, we can apply Proposition 28(ii): $\overline{M}_{\mathbb{T}_r}(f(\cdot, \theta))$ a.s. converges to $(\mu, P(f(\cdot, \theta)))$. Now, $P(f(\cdot, \theta))(x) = \mathbb{E}_x[(X_2 - \alpha_0 X_1 - \beta_0)^2 + (X_3 - \alpha_1 X_1 - \beta_1)^2] = \mathbb{E}[\varepsilon_2^2 + \varepsilon_3^2] = 2\sigma^2$. Hence, $\overline{M}_{\mathbb{T}_r}(f(\cdot, \theta))$ a.s. converges to $2\sigma^2$. Thus, it is enough to prove that $D_r$ a.s. tends to 0. Let us write $\theta = (\theta_1, \theta_2, \theta_3, \theta_4) \in \Theta = (-1, 1) \times \mathbb{R} \times (-1, 1) \times \mathbb{R}$. From the Taylor–Lagrange formula, for any $\Delta \in \mathbb{R}^3$ and $\theta, \theta' \in \Theta$, we can find $\lambda \in (0, 1)$ such that

$$f(\Delta, \theta') - f(\Delta, \theta) = \sum_{j=1}^{4} (\theta'_j - \theta_j) \partial_{\theta_j} f(\Delta, \theta + \lambda(\theta' - \theta)).$$

Now, observing that $f$ is a polynomial of global degree 4 and of degree 2 in each $\theta_j$, we can find $g \in C_{\text{pol}}(\mathbb{R}^3)$ such that for all $j \in \{1, 2, 3, 4\}$, $\Delta \in \mathbb{R}^3$ and $\theta \in \Theta$, $|\partial_{\theta_j} f(\Delta, \theta)| \leq g(\Delta)(1 + \|\theta\|)$. Therefore, for all $r \in \mathbb{N}$,

$$|D_r| \leq \frac{1}{2} \|\hat{\theta}^r - \theta\| \frac{1}{|\mathbb{T}_r|} \sum_{i \in \mathbb{T}_r} g(\Delta_i)(1 + \|\theta\| + \|\hat{\theta}^r - \theta\|)$$

$$= \frac{1}{2} \|\hat{\theta}^r - \theta\|(1 + \|\theta\| + \|\hat{\theta}^r - \theta\|) \overline{M}_{\mathbb{T}_r}(g).$$

From Proposition 28(ii), $\overline{M}_{\mathbb{T}_r}(g)$ a.s. converges. Besides, Proposition 31 states that $\|\hat{\theta}^r - \theta\|$ a.s. tends to 0. As a consequence, so does $D_r$. This completes the proof for $\sigma^2$. The proof for $\rho$ is very similar. $\square$

3.3. *Detection of cellular aging.* As explained in [13], detecting cellular aging boils down, in the BAR model (1), to rejecting hypothesis $H_0 = \{(\alpha_0, \beta_0) = (\alpha_1, \beta_1)\}$. Let us now build a statistical test that allows us to segregate between $H_0$ and its alternative $H_1 = \{(\alpha_0, \beta_0) \neq (\alpha_1, \beta_1)\}$. Wald's test is well adapted to the situation. We write $\hat{\mu}_{1,r} = \mu_1(\hat{\theta}_r)$ and $\hat{\mu}_{2,r} = \mu_2(\hat{\theta}_r, \hat{\sigma}_r)$ [recall (19)].

PROPOSITION 35. *In the BAR model* (1), *assume that the distribution of the ancestor $X_1$, $\nu$, has finite moments of all orders. Then the test statistic*

$$\chi_r^{(1)} = \frac{|\mathbb{T}_r|}{2\hat{\sigma}_r^2(1 - \hat{\rho}_r)} \{(\hat{\alpha}_0^r - \hat{\alpha}_1^r)^2 (\hat{\mu}_{2,r} - \hat{\mu}_{1,r}^2) + ((\hat{\alpha}_0^r - \hat{\alpha}_1^r)\hat{\mu}_{1,r} + \hat{\beta}_0^r - \hat{\beta}_1^r)^2\}$$

*converges in distribution to $\chi^2(2)$, the $\chi^2$ distribution with two degrees of freedom, under $H_0$, and almost surely diverges to $+\infty$ under $H_1$.*

PROOF. Recall that $\theta = (\alpha_0, \beta_0, \alpha_1, \beta_1)$. Let us set $g(\theta) = (\alpha_0 - \alpha_1, \beta_0 - \beta_1)^t$. Then $H_0 = \{g(\theta) = 0\}$. From Proposition 33, $|\mathbb{T}_r|^{1/2}(\hat{\theta}^r - \theta)$ converges in law to $\mathcal{N}_4(0, \Sigma')$ so that $|\mathbb{T}_r|^{1/2}(g(\hat{\theta}^r) - g(\theta))$ converges in law to $\mathcal{N}_2(0, \Sigma'')$, with

$$\Sigma'' = dg(\theta)\Sigma' dg(\theta)^t$$
$$= \sigma^2 (I_2 \quad -I_2) \begin{pmatrix} K & \rho K \\ \rho K & K \end{pmatrix} \begin{pmatrix} I_2 \\ -I_2 \end{pmatrix} = 2\sigma^2(1-\rho)K.$$

Under $H_0$, $g(\theta) = 0$ so that $|\mathbb{T}_r|^{1/2} g(\hat{\theta}^r)$ converges in law to $G \sim \mathcal{N}_2(0, \Sigma'')$. Now, from (18), $K = K(\theta, \sigma)$ is a continuous function of $(\theta, \sigma) \in \Theta \times \mathbb{R}_+^*$ so that, letting $\hat{K}_r = K(\hat{\theta}^r, \hat{\sigma}_r)$, Propositions 31 and 34 imply that $\hat{\Sigma}''_r \equiv 2\hat{\sigma}_r^2(1 - \hat{\rho}_r)\hat{K}_r$ converges in probability to $\Sigma''$. By continuity of $G \mapsto \Sigma''^{-1/2}G$, Slutsky's theorem shows that $|\mathbb{T}_r|^{1/2} \hat{\Sigma}''_r{}^{-1/2} g(\hat{\theta}^r)$ converges in law to $\mathcal{N}_2(0, I_2)$. In particular,

$$\||\mathbb{T}_r|^{1/2}\hat{\Sigma}''_r{}^{-1/2} g(\hat{\theta}^r)\|^2 = |\mathbb{T}_r| g(\hat{\theta}^r)^t \hat{\Sigma}''_r{}^{-1} g(\hat{\theta}^r) = \frac{|\mathbb{T}_r|}{2\hat{\sigma}_r^2(1-\hat{\rho}_r)} g(\hat{\theta}^r)^t \hat{K}_r^{-1} g(\hat{\theta}^r)$$

equals $\chi_r^{(1)}$ and converges to the $\chi^2$ distribution with two degrees of freedom.

Under $H_1$, $\chi_r/|\mathbb{T}_r| = g(\hat{\theta}^r)^t \hat{\Sigma}''_r{}^{-1} g(\hat{\theta}^r)$ a.s. converges to $g(\theta)^t \Sigma''^{-1} g(\theta) > 0$ so that $\chi_r^{(1)}$ a.s. diverges to $+\infty$. $\square$

The same technique may be used to test $\{\alpha_0 = \alpha_1\}$:

PROPOSITION 36. *In the BAR model* (1), *assume that the distribution of the ancestor $X_1$, $\nu$, has finite moments of all orders. Then the test statistic*

$$\chi_r^{(2)} = |\mathbb{T}_r|(\hat{\alpha}_0^r - \hat{\alpha}_1^r)^2 \frac{\hat{\mu}_{2,r} - \hat{\mu}_{1,r}^2}{2\hat{\sigma}_r^2(1 - \hat{\rho}_r)}$$

*converges in distribution to $\chi^2(1)$ under $\{\alpha_0 = \alpha_1\}$ and a.s. diverges to $+\infty$ under $\{\alpha_0 \neq \alpha_1\}$.*

The same can be done for testing $\{\beta_0 = \beta_1\}$. Proposition 31 provides natural statistics to test $\{\alpha_0 = 0\}$, $\{\alpha_1 = 0\}$, $\{\beta_0 = 0\}$, $\{\beta_1 = 0\}$ and $\{\beta_0 = \beta_1\}$. We do not give details for the sake of brevity.

We now present a statistical test that allows us to differentiate between $H'_0 = \{\beta_0/(1-\alpha_0) = \beta_1/(1-\alpha_1)\}$ and its alternative $H'_1$. This allows to test if the two fixed points corresponding to the two affine regressions of the BAR model (1) are equal. This may happen even if $(\alpha_0, \beta_0) \neq (\alpha_1, \beta_1)$. Rejecting $H'_0$ means accepting that the new pole and the old pole populations are even distinct in mean. Again, we use Wald's test, since $H'_0 = \{g(\theta) = 0\}$ with $g(\theta) = \beta_0/(1-\alpha_0) - \beta_1/(1-\alpha_1)$. The proof is obvious and not detailed here.

PROPOSITION 37. *In the BAR model* (1), *assume that the distribution of the ancestor $X_1$, $\nu$, has finite moments of all orders. Let*

$$dg(\theta) = \left( \frac{\beta_0}{(1-\alpha_0)^2}, \frac{1}{1-\alpha_0}, \frac{-\beta_1}{(1-\alpha_1)^2}, \frac{1}{1-\alpha_1} \right)$$

*be the gradient of $g$ and $\hat{s}_r^2 = dg(\hat{\theta}_r)\hat{\Sigma}'_r dg(\hat{\theta}_r)^t$, where $\hat{\Sigma}'_r$ is $\Sigma'$ evaluated in $(\hat{\theta}_r, \hat{\sigma}_r, \hat{\rho}_r)$, and $\Sigma'$ is defined in* (20). *Then the test statistic*

$$\chi_r^{(3)} = \frac{|\mathbb{T}_r|}{\hat{s}_r^2} \left( \frac{\hat{\beta}_0^r}{1-\hat{\alpha}_0^r} - \frac{\hat{\beta}_1^r}{1-\hat{\alpha}_1^r} \right)^2$$

*converges in distribution to $\chi^2(1)$ under $H_0'$ and a.s. diverges to $+\infty$ under $H_1'$.*

In the case when $\alpha_0 = \alpha_1 = 0$, testing $H_0$ or $H_0'$ boils down to testing $\{\beta_0 = \beta_1\}$:

PROPOSITION 38. *In the BAR model* (1), *assume that $\alpha_0 = \alpha_1 = 0$ and that the distribution of the ancestor $X_1$, $\nu$, has finite moments of all orders. Then the test statistic*

$$\xi_r = \frac{1}{\hat{\sigma}_r \sqrt{2|\mathbb{T}_r|(1-\hat{\rho}_r)}} \sum_{i \in \mathbb{T}_r} (X_{2i} - X_{2i+1})$$

*converges in distribution to $\mathcal{N}(0,1)$ under $\{\beta_0 = \beta_1\}$ and almost surely tends to $+\infty$ (resp. $-\infty$) under $\{\beta_0 > \beta_1\}$ (resp. $\{\beta_0 < \beta_1\}$).*

PROOF. Let $f(x,y,z) = y - z$. Observe that $f \in C_{\text{pol}}(\mathbb{R}^3)$ and that $Pf(x) = \beta_0 - \beta_1$, since $\alpha_0 = \alpha_1 = 0$.

Let us assume that $\beta_0 = \beta_1$. Then $Pf = 0$ and Proposition 27(iii) implies that $\hat{\sigma}_r \sqrt{2(1-\hat{\rho}_r)}\xi_r$ converges in distribution to $\mathcal{N}(0,s^2)$, where $s^2 = (\mu, Pf^2) = 2\sigma^2(1-\rho)$. Now, $(\hat{\sigma}_r^2, \hat{\rho}_r)$ a.s. converges to the constant $(\sigma^2, \rho)$ so that, with Slutsky's theorem, $\xi_r$ converges in distribution to $\mathcal{N}(0,1)$.

Let us now assume that $\beta_0 \neq \beta_1$. Proposition 28 states that the sequence $|\mathbb{T}_r|^{-1/2} \hat{\sigma}_r \sqrt{2(1-\hat{\rho}_r)}\xi_r = |\mathbb{T}_r|^{-1} \sum_{i \in \mathbb{T}_r} f(\Delta_i)$ a.s. converges to $(\mu, Pf) = \beta_0 - \beta_1 \neq 0$. Since $\hat{\sigma}_r \sqrt{2(1-\hat{\rho}_r)}$ converges in probability to $\sigma\sqrt{2(1-\rho)} > 0$, we conclude that $\xi_r$ a.s. diverges as $|\mathbb{T}_r|^{1/2}$, to $+\infty$ if $\beta_0 > \beta_1$, and to $-\infty$ if $\beta_0 < \beta_1$. □

REMARK 39. Note that in the model where $\alpha_0 = \alpha_1 = 0$, $\chi_r^{(1)}$ would read

$$\chi_r^{(1)} = \frac{|\mathbb{T}_r|}{2\hat{\sigma}_r^2(1-\hat{\rho}_r)} (\hat{\beta}_0^r - \hat{\beta}_1^r)^2$$

and is thus equal to $\xi_r^2$. The latter test looks like the ones Stewart et al. performed in [22]: it focuses on the differences $X_{2i} - X_{2i+1}$ between sisters. But it is relevant only in the case when the correlation parameters $\alpha_0$ and $\alpha_1$ are zero, that is, in a dynamics with no memory. Now, the data analysis strongly rejects this assumption, as we shall see in the next section.

REMARK 40. Note that all the propositions and lemmas of Section 3 remain true if $((\varepsilon_{2i}, \varepsilon_{2i+1}), i \geq 1)$ is only assumed to be a sequence of (non-necessarily Gaussian) i.i.d. bivariate random vectors with finite moments of all orders, that is, $\mathbb{E}[\varepsilon_{2i}^{2m} + \varepsilon_{2i+1}^{2m}] < \infty$ for all $m \in \mathbb{N}$, with covariance matrix given by (2). In such a general case, $(\hat{\theta}_r, \hat{\sigma}_r^2, \hat{\rho}_r)$ has no reason to be the maximum likelihood estimator of $(\theta, \sigma, \rho)$.

3.4. *Data numerical analysis.* We now perform the estimation and test procedures on Stewart et al.'s data (see Guyon [12] for more detailed results). The data consists of 95 films, and each film should be seen as an incomplete binary tree of growth rates. How do we compute the estimators and test statistics? According to the above presentation, we should restrict the observation to the bigger complete subtree $\mathbb{T}_{r+1}$. We actually take into account all the observations, noting that:

- very few cells are observed in a generation, say, $r$, when generation $r-1$ is not completely observed,
- cells observed in the last generation are assumed to be the result of a random permutation $\Pi$, independent of $X$; this should be correct as a first approximation.

Figure 5 gives the global empirical distribution of the residues $\hat{\varepsilon}$ over the 95 films. We have separated new poles' residues (left) from old poles' ones (right), and fitted to normal distributions. Both histograms are close to Gaussian laws.

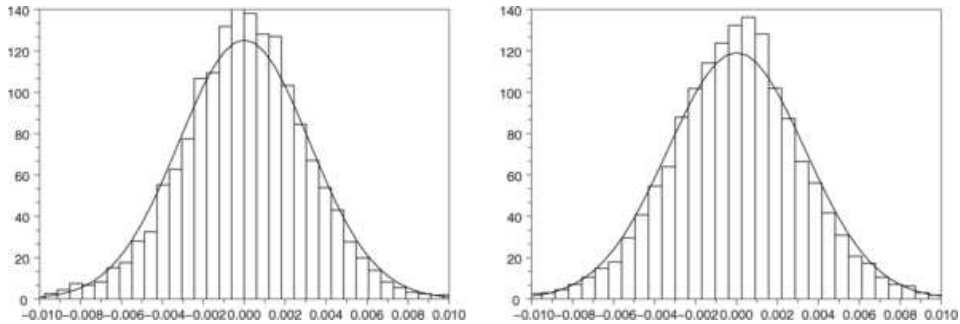

FIG. 5. *Histograms of the residues $\hat{\varepsilon}_{2n}$ (new poles, left) and $\hat{\varepsilon}_{2n+1}$ (old poles, right), and their fit to Gaussian distributions.*

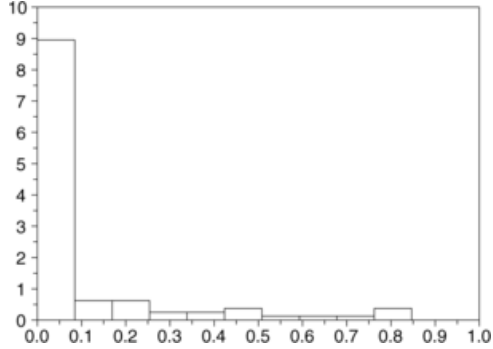

FIG. 6. $H_0 = \{(\alpha_0, \beta_0) = (\alpha_1, \beta_1)\}$. Histogram of the p-values $\mathbb{P}(\chi^2(2) \geq \chi^{(1)}_{\mathrm{obs}})$.

Figure 6 shows that $H_0$ can be strongly rejected. This indicates that the dynamics of the growth rate of the old pole offspring is different from that of the new pole offspring. The nullity of any parameter ($\alpha_0$, $\beta_0$, $\alpha_1$ or $\beta_1$) can be strongly rejected as well. This enlightens the relevance of a Markovian modelization with memory one: the mother cell is a significant predictor of offspring growth rate in general.

Besides, we cannot reject the hypothesis that both $\alpha$'s are equal on the one hand, and that both $\beta$'s are equal on the other hand. But we strongly reject that both fixed points, namely, $\gamma_0 = \beta_0/(1-\alpha_0)$ and $\gamma_1 = \beta_1/(1-\alpha_1)$, are equal; see Figure 7. Hence, the parametrization $(\alpha, \gamma)$, which makes more physical sense than the parametrization $(\alpha, \beta)$, has the following advantage: with no assumption on the $\alpha$'s, we can detect aging by looking only at the $\gamma$'s, which we cannot do with the $\beta$'s. It also means that the new poles and the old poles are not only different in distribution, but also in mean.

The scatter plot in Figure 7 indicates that $\gamma_0 > \gamma_1$. More precisely, the line $\gamma_0 = \gamma_1 + \delta$ fits well the data with $\delta$ significantly positive. Numerically, $\delta \approx 0.0012 \pm 0.0011$, or $0.0011 \pm 0.0008$ if we delete the two aberrant points in Figure 7 (right). This may be seen as statistical evidence of aging in *E. Coli*, since, on average, old pole cells grow slower than the new pole cells, which is characteristic of aged individuals. Quantitatively, they seem to grow 3% slower (we may speak in terms of percentage since the range of values of $\gamma$'s is narrow). This result is close to Stewart et al.'s original calculations, since in [22] they estimated this ratio to be around 2%.

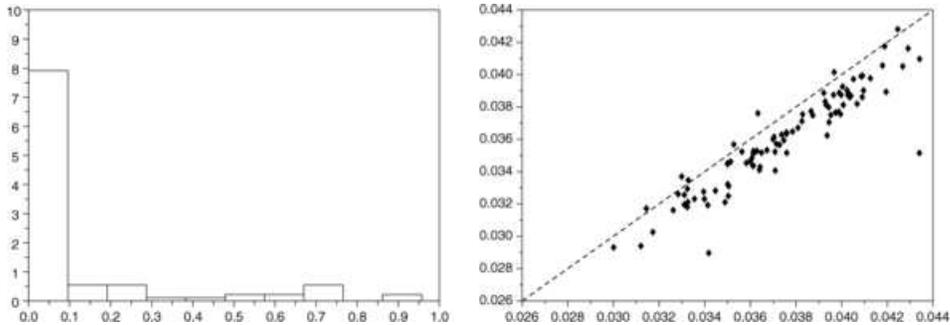

FIG. 7. $H'_0 = \{\beta_0/(1-\alpha_0) = \beta_1/(1-\alpha_1)\}$; left: histogram of the p-values $\mathbb{P}(\chi^2(1) \geq \chi^{(3)}_{\mathrm{obs}})$; right: $\hat{\beta}_0/(1-\hat{\alpha}_0)$ on the x-axis, $\hat{\beta}_1/(1-\hat{\alpha}_1)$ on the y-axis; the dashed line is the diagonal.

**Acknowledgments.** I would like to thank J.-F. Delmas for discussions in which he suggested many ideas and for thorough readings of the paper. I also wish to thank F. Taddéi, E. J. Stewart, A. Lindner and G. Paul for a pleasant and efficient collaboration, and for providing us with their original data.

CERMICS-ENPC
6 ET 8 AVENUE BLAISE PASCAL
CITÉ DECARTES CHAMPS-SUR-MARNE
77455 MARNE LA VALLÉE CEDEX 2
FRANCE
E-MAIL: julien.guyon@cermics.enpc.fr